\documentclass[a4paper,reqno,12pt]{article}
\usepackage{shuffle}
\usepackage{hyperref}
\usepackage{amsmath,amsfonts,amssymb}
\usepackage{fancyhdr}
\usepackage{mathrsfs,amssymb} 
\usepackage{amsmath}
\usepackage{amssymb}
\usepackage{empheq}
\usepackage{enumerate}
\usepackage{tikz}
\usepackage{float}
\usepackage[a4paper]{geometry}
\usepackage[latin1]{inputenc}
\usepackage[T1]{fontenc}
\usepackage{amsthm}
\usepackage{dsfont}
\usepackage[nottoc]{tocbibind}
\usepackage{color}
\usepackage{fancyhdr,blindtext}
\usepackage[arrow, matrix, curve]{xy}
\usepackage{array}
\newcolumntype{x}[1]{>{\centering\arraybackslash\hspace{0pt}}m{#1}}

\fancypagestyle{plain}{%
\fancyhf{}}

\renewcommand{\sectionmark}[1]%
{\markboth{#1}{}}

\makeindex

\theoremstyle{definition}
\newtheorem{ex}{\bfseries \upshape Example}[section]
\newtheorem{dfn}[ex]{\bfseries \upshape Definition}
\newtheorem{rem}[ex]{\bfseries \upshape Remark}
\newtheorem{conj}[ex]{\bfseries \upshape Conjecture}
\newtheorem{prop}[ex]{\bfseries \upshape Proposition}
\newtheorem{lem}[ex]{\bfseries \upshape Lemma}
\newtheorem{thm}[ex]{\bfseries \upshape Theorem}
\theoremstyle{plain}

\newtheorem{cor}[ex]{\bfseries \upshape Corollary}

\newenvironment{prf}{\begin{proof}[{\bf Proof}]}{\end{proof}}

\newcommand{\N}{\ensuremath{\mathds{N}}}	
\newcommand{\Z}{\ensuremath{\mathds{Z}}}	
\newcommand{\Q}{\ensuremath{\mathds{Q}}}	
\newcommand{\R}{\ensuremath{\mathbb{R}}}	
\newcommand{\C}{\ensuremath{\mathds{C}}}

\newcommand{\Li}{\operatorname{Li}}
\newcommand{\Lit}{\widetilde{\operatorname{Li}}}
\newcommand{\dif}{ \operatorname{d}}
\newcommand{\gr}{ \operatorname{gr}}
\newcommand{\grw}{ \operatorname{gr}^{\operatorname{W}}}
\newcommand{\grl}{ \operatorname{gr}^{\operatorname{L}}}
\newcommand{\grwl}{ \operatorname{gr}^{\operatorname{W},\operatorname{L}}}
\newcommand{\fil}{ \operatorname{Fil}}
\newcommand{\filw}{ \operatorname{Fil}^{\operatorname{W}}}
\newcommand{\fille}{ \operatorname{Fil}^{\operatorname{L}} }
\newcommand{\filwle}{ \operatorname{Fil}^{\operatorname{W},\operatorname{L}}}

\newcommand{\qe}[1]{\frac{q^{#1}}{1-q^{#1}}}

\DeclareMathOperator{\Sl}{SL}

\DeclareMathOperator{\MD}{\mathcal{MD}}
\DeclareMathOperator{\MDA}{q\mathcal{M}\mathcal{Z}}
\DeclareMathOperator{\MDE}{\mathcal{MD}^\textrm{even}}
\DeclareMathOperator{\MQ}{\mathcal{Q}}

\DeclareMathOperator{\MZ}{\mathcal{MZ}}

\newcommand{\thmref}[1] {Theorem \ref{#1}}
\newcommand{\lemref}[1] {Lemma \ref{#1}}
\newcommand{\propref}[1] {Proposition \ref{#1}}
\newcommand{\remref}[1] {Remark \ref{#1}}

\newcommand{\tabref}[1] {Table \ref{#1}}
\newcommand{\exref}[1] {Example \ref{#1}}
\newcommand{\defref}[1] {Definition \ref{#1}}

\newcommand{\blue}[1]{{\color{blue}{#1}}}
\numberwithin{equation}{section}

\begin{document}
\title{{ \bf The algebra of generating functions for multiple divisor sums and applications to multiple zeta values}}
\author{{\sc Henrik Bachmann}, {\sc Ulf K\"uhn}}
\date{\today}
\maketitle
\textheight215mm

\begin{abstract}
We study the algebra $\MD$ of generating function for multiple divisor sums and its connections to multiple zeta values. The generating functions for multiple divisor sums are formal power series in $q$ with coefficients in $\Q$ arising from the calculation of the Fourier expansion of multiple Eisenstein series.  We show that the algebra $\MD$ is a   filtered algebra equipped with a derivation and use this derivation to prove linear relations in $\MD$.  The (quasi-)modular forms for the full modular group $\Sl_2(\Z)$ constitute a sub-algebra of $\MD$ this also yields linear relations in $\MD$. Generating functions of multiple divisor sums can be seen as a $q$-analogue of multiple zeta values. Studying a certain map from this algebra into the real numbers we will derive a new explanation for relations between multiple zeta values, including those in length $2$, coming from modular forms. 
\end{abstract}

\tableofcontents

\section{Introduction}
Multiple zeta values are natural generalizations of the Riemann zeta values that are defined for integers $s_1 > 1$ and $s_i \geq 1$ for $i>1$ by
\[ \zeta( s_1 , \dots , s_l ) :=  \sum_{n_1 >n_2 >\dots>n_l>0 } \frac{1}{n_1^{s_1} \dots n_l^{s_l} } \,. \]
Because of its occurence in various fields of mathematics and physics these real numbers 
are of particular interest.
The $\Q$-vector space of all multiple zeta values of weight $k$ is then given by 
\[ \MZ_k :=\big <\,\zeta(s_1,\dots,s_l) \, \big| \, s_1 + \dots + s_l = k 
\textrm{ and } l>0 \big>_\Q.  \]
It is well known that  the product of two multiple zeta values can be written as a linear combination of multiple zeta values of the same weight by using the stuffle or shuffle relations. Thus they generate a $\Q$-algebra $\MZ$. There are  beautiful conjectures about the dimensions of finite dimensional subspaces 
of $\MZ$ determined by the weight and the depth filtration.

In \cite{gkz} Gangl, Kaneko and Zagier introduced double Eisenstein series, which were generalized to
multiple Eisenstein series in \cite{hb}. These series are sums over certain positive sectors in the multiple product of a lattice. They give natural generalizations of the well-known Eisenstein series from the theory of modular forms similar as the multiple zeta values generalize special values of the Riemann zeta function. These functions do by construction satisfy the stuffle relations. But 
due to convergence problems the shuffle relation 
needs some modification; it seems to hold up to an error term which involves derivatives. 
The motivation behind this article is the idea to understand these corrections algebraically, although this will not be discussed here furthermore (c.f. \cite{BBK}, \cite{BT}).
It has been  shown in \cite{hb} that multiple Eisenstein series have a Fourier expansion, which decomposes
as a $\MZ$-linear combination of generating functions for multiple divisor sums $[s_1,\dots,s_l]$ which we also refer to as brackets in this paper.  For example the double Eisenstein series $G_{4,4}$  and the triple Eisenstein series $G_{3,2,2}$ are given by
 \begin{align*}
G_{4,4}(\tau) =& \zeta(4,4) + 20 \zeta(6) (2 \pi i)^2 [2](q_\tau) + 
3 \zeta(4) (2 \pi i)^4 [4](q_\tau) + (2\pi i)^8 [4,4](q_\tau) \,,\\
G_{3,2,2}(\tau) =& \zeta(3,2,2) + \left( \frac{54}{5} \zeta(2,3) + \frac{51}{5} \zeta(3,2) \right) (2 \pi i)^2 [2](q_\tau) + \frac{16}{3} \zeta  (2,2) (2 \pi  i)^3 [3](q_\tau) \\
&+3 \zeta(3) (2\pi i)^4 [2,2](q_\tau) + 4 \zeta(2) (2 \pi i)^5 [3,2](q_\tau) + (2 \pi i)^7[3,2,2](q_\tau)\,,
\end{align*}
where
$\tau \in \mathbb{H}$, $q_\tau= \exp(2 \pi i \tau)$ and
the brackets $[s_1,\dots,s_l]$ are kind of a combinatorial object\footnote{In \cite{gkz} certain linear combinations of these functions were called combinatorial Eisenstein series} that will be described now. 
As a generalization of the classical divisor sums we define for natural numbers $r_1,\dots,r_l \in \N_0 =\{0,1,2,\dots\}$  the \emph{multiple divisor sum} by
\begin{equation} \label{def:sigma} \sigma_{r_1,\dots,r_l}(n) = \sum_{\substack{u_1 v_1 + \dots + u_l v_l = n\\u_1 > \dots > u_l >0}} v_1^{r_1} \dots v_l^{r_l} \,. \end{equation}

  For any integers $s_1,\dots,s_l>0$ the generating function for  the multiple divisor sum $\sigma_{s_1-1,\dots,s_l-1}$ is defined by the formal power series
\[ [s_1,\dots,s_l] := \frac{1}{(s_1-1)! \dots (s_l-1)!} \sum_{n>0} \sigma_{s_1-1,\dots,s_l-1}(n) q^n \,\,\in \Q[[q]] \,. \]
Here and in the following, we will simply write $[s_1,\dots,s_l]$ instead of 
$[s_1,\dots,s_l](q)$.
We refer to these generating functions of multiple divisor sums also as \emph{brackets}.\footnote{
The brackets $[2,\dots,2]$ were in the context of partitions already studied by P.A.~MacMahon (see \cite{mm}) and named generalized divisor sums.  It was shown in \cite{ar} that these are quasi-modular forms, see also \remref{rem:md-partition}}
. 

\begin{ex} We give a few examples:
\begin{align*}
[2] &= q + 3q^2 + 4q^3 + 7q^4 + 6q^5 + 12q^6 + 8q^7 + 15q^8 + \dots \,, \\
 [4,2] &=\frac{1}{6} \left( q^3 + 3q^4 + 15q^5 + 27q^6 + 78q^7 + 135q^8 +  \dots  \right)\,,  \\
 [4,4,4] &= \frac{1}{216} \left( q^6 + 9q^7 + 45q^8 + 190q^9 + 642q^{10} + 1899q^{11}  +\dots \right)\,,  \\
 [3,1,3,1] &= \frac{1}{4}\left( q^{10} + 2q^{11} + 8q^{12} + 16q^{13} + 43q^{14} + 70q^{15}+\dots \right) \,, \\
 [1,2,3,4,5]&=\frac{1}{288} \left( q^{15} + 17q^{16} + 107q^{17} + 512q^{18} + 1985q^{19} + \dots \right)\,.
\end{align*}
\end{ex}
Notice that the first non vanishing coefficient of $q^n$ in $[s_1,\dots,s_l]$ appears at \linebreak  $n = \frac{l (l+1)}{2}$, because it belongs to the "smallest" possible partition 
\[l \cdot 1 + (l-1) \cdot 1 + \dots + 1\cdot 1 = n \,, \]
i.e. $u_j= j$ and $v_j = 1$ for $1 \leq j \leq l$.  
The number $k=s_1+\dots+s_l$ is called the \emph{weight} of  $[s_1,\dots,s_l]$ and $l$ denotes the \emph{length}. These numbers satisfy $l \le k$. 

\begin{dfn}
We define the vector space $\MD$ to be the $\Q$ vector space generated by $[\emptyset]=1 \in \Q[[q]]$ and all 
brackets $[s_1,\dots,s_l]$. 
On $\MD$ we have the increasing filtration $\filw_{\bullet}$
given by the weight and the  increasing filtration $\fille_{\bullet}$
given by the length, i.e., we have
\begin{align*}
\filw_k(\MD) &:=  \big<[s_1,\dots,s_l] \,\big|\, s_1+\dots+s_l \le k \,\big>_{\Q}\\
\fille_l(\MD) &:=  \big<[s_1,\dots,s_r] \,\big|\, r\le l \,\big>_{\Q}\,.
\end{align*}
If we consider the length and weight filtration at the same time we use the short notation $\filwle_{k,l} := \filw_k \fille_l$.
As usual we set
\begin{align*}
\grw_k(\MD) &:=   \filw_k(\MD) \slash \filw_{k-1}(\MD) \\
\grl_l(\MD) &:=  \fille_l(\MD) \slash \fille_{l-1}(\MD) \,.
\end{align*}
and as above $\grwl_{k,l} := \grw_k \grl_l$. 
 
\end{dfn}


For example for even $k\geq4$ the Eisenstein series $G_{k}$, which are well-known to be modular forms of weight $k$ for the group $\Sl_2(\Z)$, are elements in this vector spaces, because they satisfy
\[ G_{k} =  \frac{\zeta(k)}{(-2 \pi i)^k} + \frac{1}{(k-1)!} \sum_{n>0} \sigma_{k-1}(n) q^n  =    -\frac{1}{2} \frac{B_{k}}{k!}[ \emptyset ]  + [k] 
\in \filw_k(\MD),\] 
also the quasi-modular form $G_2$ of weight $2$ is an element of $\filw_2(\MD)$.
Our first result is 
 
\begin{thm} \label{thm:md-algebra} 
\emph{The 
$\Q$-vector space $\MD$
has the structure of a bifiltered $\Q$-Algebra $(\MD,\, \cdot,\,\filw_{\bullet},\,\fille_{\bullet})$,   where the multiplication is the natural multiplication of formal power series and the filtrations $\filw_{\bullet}$ and  $\fille_{\bullet}$ are induced by the $weight$ and $length$,  in particular  
\[ 
\filwle_{k_1,l_1}(\MD) \cdot  \filwle_{k_2,l_2}(\MD) \subset \filwle_{k_1+k_2,l_1+l_2}(\MD). 
\label{prop:alg}\]}
\end{thm}

\begin{rem} In fact we prove that this product on $\MD$ is a   quasi-shuffle product in the sense of 
Hofmann and Ihara \cite{hi}.
\end{rem}


\begin{ex}
The first products of brackets are given by
\begin{align}
 [1] \cdot [1] &= 2[1,1] + [2] - [1] \,, \label{example:prod11} \\
[1] \cdot [2] &= [1,2]+[2,1]+[3] -\frac{1}{2}[2] \,,  \label{example:prod12} \\
[1] \cdot [2,1] &= [1,2,1]+2 [2,1,1] - \frac{3}{2} [2,1]+[2,2]+[3,1] \,.\label{example:prod121}
 \end{align}
\end{ex}

For small weight $k$ or at least a small $l$ length we can compute a sufficiently large number 
of the Fourier coefficients of a bracket. We can therefore determine
lower bounds for the number of linearly independent elements in $\filwle_{k,l}(\MD)$, in order to do so we need to check that the matrix of with rows given by the Fourier coefficients of each element has a sufficient high rank. 

\begin{thm} \label{thm:lowdim_fil_md}  
We have the following exact values or
lower bounds for $\dim_\Q \fil^{\operatorname{W},\operatorname{L}}_{k,l} (\MD)  $
\begin{table}[H]\footnotesize
\begin{center}
\begin{tabular} { c | c | c |c|c|c|c|c|c|c|c|c|c|c|c}
$k \backslash l$&0&1&2&3&4&5&6&7&8&9&10&11 \\ \hline
0& \bf{1}	&         &         & 			  & 		& 		& 		&		  & 		& 		& &		\\ \hline
1& \bf{1}	& \bf{2}  &         & 			  & 		& 		& 		&  		& 		& 		& 	&	\\ \hline
2& \bf{1}	& \bf{3}  & \bf{4}  & 			  & 		& 		& 		&  		&  		&			& &		\\ \hline
4& \bf{1}	& \bf{4}  & \bf{7}  & \bf{8}  & 		& 		& 		& 		&  		&			& &		\\ \hline
3& \bf{1}	& \bf{5}  & \bf{10} & \bf{14} & \bf{15}	& 		&			& 		&  		& 		& &		\\ \hline
5& \bf{1}	& \bf{6}  & \bf{14} & \bf{22} & \bf{27}	& \bf{28}	&			&			&  		&			& &		\\ \hline
6& \bf{1}	& \bf{7}  & \bf{18} & \bf{32}			& \blue{44}	& \blue{50}	& \blue{51}	&			& 		&			& 	&	\\ \hline
7& \bf{1}	& \bf{8}  & \bf{23} & \blue{ 44}			& \blue{ 67}	& \blue{84}	& \blue{91}	& \blue{92}	& 		&			& &		\\ \hline
8& \bf{1}	& \bf{9}  & \bf{28} & \blue{ 59}			& \blue{97}	&	\blue{133}	& \blue{156}	& \blue{164}	& \blue{165}	&			& &		\\ \hline
9& \bf{1}	& \bf{10} & \bf{34} & \blue{ 76}			& \blue{135}	& \blue{200}	& \blue{254}	& \blue{284}	& \blue{293}	&	\blue{294} &		&	\\ \hline		
10& \bf{1}& \bf{11} & \bf{40} & \blue{ 97} 			& \blue{183} & \blue{290} & \blue{396} & \blue{474} & \blue{512} & \blue{522} & \blue{523}& \\ \hline
11&\bf{1}&\bf{12}&\bf{47}&\blue{120}&\blue{242}&\blue{408}&\blue{594}&\blue{760}&\blue{869}&\blue{916}&\blue{927}&\blue{928}\\ \hline
12&\bf{1}&\bf{13}&\bf{54}&\blue{147}&\blue{313}&\blue{559}&?&?&?&?&?&?\\ \hline
13&\bf{1}&\bf{14}&\bf{62}&\blue{177}&\blue{398}&?&?&?&?&?&?&?\\ \hline
14&\bf{ 1}&\bf{15}&\bf{70}&\blue{212}&\blue{498}&?&?&?&?&?&?&?\\ \hline
15& \bf{1}&\bf{16}&\bf{79}&\blue{249}&?&?&?&?&?&?&?&?\\ \hline
\end{tabular} 
\caption{$\dim_\Q  \filwle_{k,l}(\MD)$: {\bf exact value}, \blue{lower bound}} \label{tab:filwl-md}
\end{center}
\end{table}
 
\end{thm}
The number of generators  of $\filwle_{k,l}(\MD)$ is easily calculated, thus
giving an upper bound for the dimension of this space is equivalent to give
a lower bound for the number of relations in the generators of $\filwle_{k,l}(\MD)$.
The equalities come from the fact that  we know enough relations in the cases marked
black in \tabref{tab:filwl-md}.

 For the multiple zeta values   conjecturally all linear relations are due the fact that the shuffle
and the stuffle relations give two different description of the product of two multiple zeta values,
albeit in practice there are different methods to prove distinct relations 
like the cyclic sum identity \cite{HofOhn} or  the Zagier-Ohno relation \cite{ZagOhn}.   
So far we know only one way to write a product of two brackets as a linear combination in $\MD$   and this doesn't suffice to give linear relations between elements in $\MD$. However, 
 as we will see  now,  $\MD$ has the additional structure of a differential algebra and moreover
 there are several ways to express the derivative of a bracket. 
 By now linear relations
 in $\MD$ are proved either by using derivatives and or the theory of quasi-modular forms.   
 
\begin{thm} \label{thm:derivative} The operator $\dif = q \frac{d}{dq}$ is a derivation on $\MD$, it maps $\filwle_{k,l}(\MD)$ to $\filwle_{k+2,l+1}(\MD)$.
\end{thm}

Our proof actually  allows us to derive explicit formulas for  $\dif[s_1]$ and $\dif[s_1,s_2]$. 

\begin{rem} Our formula for $\dif[k]$ may be seen as the Euler decomposition 
formula for $\MD$, since for we prove in
  \propref{cor:formularfordk}   that for $s_1+s_2=k+2$
\[
[s_1]\cdot [s_2] =
   \sum_{a+b=k+2} \left( \binom{a-1}{s_1-1}+\binom{a-1}{s_2-1} \right) [a,b] 
 - \binom{k}{s_1-1} [k+1]  + \binom{k}{s_1-1} \frac{d[k]}{k}   
   \,.\]
Frankly speaking the derivative $\dif[k]$ measures the failure of the shuffle relation for the
product of two length one bracket. 
\end{rem}

We will show now how to derive  from these formulas non trivial linear relations.

\begin{ex}(Relations from derivatives) \label{ex:relderiva}The first derivatives are given by
 \begin{align} \label{example:derivative1}
\dif [1] &= [3] + \frac{1}{2}[2] - [2,1] \,,  \\
\dif[2] &= [4] + 2 [3] - \frac{1}{6} [2] - 4 [3,1] \,,\label{example:derivative2}\\
\dif[2] &= 2[4] + [3] +\frac{1}{6} [2] -2[2,2] - 2[3,1] \,,\label{example:derivative3}\\
\dif [1,1] &=[3,1] + \frac{3}{2}[2,1] + \frac{1}{2}[1,2]+[1,3]- 2 [2,1,1]-[1,2,1] \,. \label{example:derivative4}
\end{align}
The difference of \eqref{example:derivative2} and \eqref{example:derivative3} leads to the first linear relation in
$\filw_{4}(\MD)$: 
 \begin{align}\label{eq:rel4}
  [4] =  2[2,2] - 2 [3,1] + [3]- \frac{1}{3} [2] \,.
  \end{align}
  \end{ex}
  \begin{ex}(Leibniz rule)  Since $\dif$ is a derivation it satisfies the
  Leibniz rule, e.g.,  because of  \eqref{example:prod11}
  \[  \dif[1] \cdot [1] + [1] \cdot \dif[1]
  = \dif( [1] \cdot [1]) =  \dif( 2[1,1] + [2] - [1] )  \,.  \]
 Now using \eqref{example:derivative1},  \eqref{example:derivative2} and \eqref{example:derivative4} together with the explicit description of the various products  we could alternatively prove the relation \eqref{eq:rel4}.
\end{ex}
\begin{ex}(Relations from modular forms)\label{ex:relmod}
It is a well-known fact from the theory of modular forms that $G_4^2 = \frac{7}{6} G_8$ because the space of weight $8$ modular forms for $\Sl_2(\Z)$ is one dimensional. We therefore have   
\begin{align*}
   \frac{1}{720} [4]  + [4] \cdot [4] 
&= \frac{7}{6} [8]\,.
\end{align*}
Using the product as described in \propref{prop:l2-expli} we get
\begin{align*}
 [4] \cdot [4] 
&= 2 [4,4] + [8] + \frac{1}{360} [4] - \frac{1}{1512} [2] \,,
\end{align*}
which then gives the following relation in $\filw_{8}(\MD)$:
\begin{align}\label{eq:relwt8}
[8] = \frac{1}{40} [4] - \frac{1}{252} [2] + 12 [4,4]\,.  
\end{align}
\end{ex}

Beside the methods mentioned in \exref{ex:relderiva} and \ref{ex:relmod} 
other obvious ways to get relations in weight $k$ are either to multiply a relation in weight $l$ 
by a bracket of weight $k-l$ or to take the derivative of a relation in weight $k-2$.

\begin{ex}\label{ex:oldrelations}(Relations from known relations)
If we multiply the relation \eqref{eq:rel4} in weight $4$  with $[2]$, then we obtain in
$\filw_{6}(\MD)$:
\begin{small} 
\begin{align} 
\begin{split}
\label{eq:rel6len3_1}
[6]& =  \frac{1}{20} [2] - \frac{1}{12} [3]- \frac{1}{4}[4]+[5]- \frac{4}{3} [2,2] + \frac{1}{6} [3,1]+[2,3]+2 [3,2]\\
  &+ 6 [2,2,2] -2 [3,1,2]-2 [2,3,1]-2 [3,2,1]+[2,4]-2 [3,3]+[4,2]-2 [5,1] \,. 
  \end{split}
 \end{align}
 \end{small}
If we apply $\dif$ to 
the relation \eqref{eq:rel4} in weight $4$, then we obtain in
$\filw_{6}(\MD)$:
\begin{small} 
\begin{align} \label{eq:rel6len3_2}
[6]& =  \frac{1}{20} [2]- \frac{3}{4} [3] + \frac{11}{4} [4] - 3 [5] - \frac{2}{3} [2,2]+ \frac{3}{2} [3,1] + 4 [2,3]+2 [2,4] \\
&+5 [3,2]-18 [4,1]+5 [4,2]+6 [5,1]-8 [2,3,1]-8 [3,1,2]-2 [3,2,1]+18 [4,1,1]\,. \notag
 \end{align}
 \end{small}
\end{ex}

In order to study the linear relations in the generators of $\MD$ systematically 
it is better
first to understand some of the algebra structure of $\MD$. For this purpose 
we call a brackets $[s_1,\dots,s_l]$ admissible, if $s_1>1$.     We show that 
the vector space $\MDA$ of admissible brackets is a sub algebra of $\MD$. 
In addition we prove that
$\MD$ is a polynomial ring over $\MDA$ with indeterminate $[1]$, i.e. we have 
  \[\MD = \MDA[ [1] ]\]
(see  \thmref{thm:polyadalg}). 
With this structure in our hands it is easy see that it suffices to study the 
linear relations in the generators of the quotient spaces
$\grwl_{k,l}(\MDA)$ in order to get upper bounds on the dimensions
 of all the graded or filtrated pieces of $\MDA$ 
or $\MD$.  In \thmref{thm-dklmda} we present our results in this direction.
We like to emphasize that the focus of this article is not to
 give the best possible results on the number of relations.
We expect that with a more detailed study of the kind of relations 
we can obtain so far we could derive much better results and we plan 
to come back to this in future \cite{Ba:PhDthesis}.

The notation $\MDA$ shall emphasize the relation  to q-analogues of multiple zeta values, which will be explained
now. Our algebra $\MDA$ 
is related, but not isomorphic, to a recent modification of  
multiple q zeta values as proposed in \cite{OT} or \cite{Tak-harm}, see also \remref{rem:q-vergleich}.

Define for $k\geq 0$ the map $Z_k$ on $\filw_{k}(\MDA)$ by
\begin{align*}
Z_k[s_1,\dots,s_l] = \lim_{q \to 1} (1-q)^{k} [ s_1,\dots,s_l].
\end{align*} 
We will show that with this definition we have 
\[ Z_k\left( [s_1, \dots , s_l ] \right) = \left\{
\begin{array}{cl} \zeta(s_1,\dots,s_l)\,, & k = s_1+\dots+s_l ,  \\ 0\,, &k > s_1+\dots+s_l 
\,. \end{array}   \right. 
\]

Since $\MD = \MDA[[1]]$ we can define a map $Z^{alg}_k:  \filw_{k}(\MD) \rightarrow \R[T]$ by 
\[ Z^{alg}_k\left( \sum_{j=0}^{k} g_j [1]^{k-j} \right) = \sum_{j=0}^{k} Z_j(g_j) T^{k-j} \in \R[T]\,   \]
where $g_j \in \filw_{j}(\MDA)$.  For our next result an analytical interpretation of $Z^{alg}_k$ in a broader context is the key fact.

\begin{thm}\label{thm_Zk-kernel} For the kernel of $Z^{alg}_k \in \filw_{k}(\MD)$ we have 
\begin{enumerate}[i)]
\item If for $[s_1,\dots,s_l]$ it holds $s_1+\dots+s_l<k$, then 
$Z_k^{alg}[s_1,\dots, s_l]=0$.
\item For any $f \in \filw_{k-2}(\MD) $ we have $Z_k^{alg}\dif(f)=0$, i.e.,  $\dif \filw_{k-2}(\MD)\subseteq \ker Z_k$.
\item If $f \in \filw_{k}(\MD)$ is a cusp form for $\Sl_2(\Z)$, then $Z_k^{alg}(f)=0$, i.e.
$ S_k(\Sl_2(\Z))\subseteq \ker Z_k$.
\end{enumerate}
\end{thm} 

Using 
\thmref{thm_Zk-kernel} we get as immediate consequences 
and without any difficulties the following well-known identities for multiple zeta values.
 
\begin{ex} \begin{enumerate}[i)]
\item If we apply $Z_3$ to \eqref{example:derivative1} we deduce
$\zeta(3)=\zeta(2,1)$.
\item If we apply  $Z_4$ to \eqref{example:derivative2} and \eqref{example:derivative3} 
we deduce $\zeta(4)= 4 \zeta(3,1) = \frac{4}{3} \zeta(2,2)$.
\item  The identity \eqref{example:derivative4} reads
 in $\MDA[[1]]$ as
\[ \dif[1,1]= \left([3]-[2,1]+\frac{1}{2}[2]\right)\cdot[1]  
+ 2 [3,1] - \frac{1}{2}[4] - \frac{1}{2}[2,1] - \frac{1}{2}[3] + \frac{1}{3}[2]\,.\,\]
Applying $Z^{alg}_4$ we deduce again the two relations
$\zeta(3)=\zeta(2,1)$ and $4\zeta(3,1)=\zeta(4)$, since
by \thmref{thm_Zk-kernel} we  have 
\[ Z^{alg}_4(\dif[1,1])=\left( \zeta(3) - \zeta(2,1) \right) T - \frac{1}{2}\zeta(4) + 2 \zeta(3,1) = 0 \,. \] 
\item If we apply $Z_8$ to \eqref{eq:relwt8} we deduce $\zeta(8)= 12 \zeta(4,4)$.
\item As an application of \thmref{thm:lowdim_fil_md}  we can prove 
for the cusp form $\Delta \in S_{12}(\Sl_2(\Z))$   the representation 
\begin{align}\label{eq:deltal2}
 \frac{1}{2^6\cdot 5 \cdot 691}  \Delta  &=  168 [5,7]+150 [7,5]+28 [9,3] \notag \\
&+\frac{1}{1408} [2] - \frac{83}{14400}[4] +\frac{187}{6048} [6] - \frac{7}{120} [8] - \frac{5197}{691} [12] \,.
\end{align}
Letting $Z_{12}$ act on both sides of \eqref{eq:deltal2} one obtains the relation
\[ \frac{5197}{691} \zeta(12) =  168 \zeta(5,7)+150 \zeta(7,5) + 28 \zeta(9,3) \,. \] 
\end{enumerate}
\end{ex}
Finally we point to the fact that the last identity coming from the cusp form $\Delta$
has been obtained via  period polynomials in \cite{gkz}.  A remarkable fact of this
relation is that   
it is not provable within the double shuffle relations in weight $12$ and depth $2$ 
alone, since also the extended double shuffle relations are needed for its proof.

This article contains results that will be part of  the dissertation project  by the first author.

We thank O.~Bouillot, F.~Brown, J.~Burgos, H.~Gangl, O.~Schnetz, D.~Zagier, J.~Zhao  and W.~Zudilin 
for their interest in our work and for helpful remarks.

\section{The algebra of generating function of multiple divisor sums}

The proof of \thmref{thm:md-algebra} will occupy  this section. First we 
consider products of polylogarithms at negative integers. This will give us an explicit formula for the product of two brackets. 

\begin{rem} \label{rem:md-partition} We start with a remark on where brackets also
have appeared before.
In the following we will write $\{a\}^l$ for a length $l$ sequence $a,\dots,a$. 

\begin{enumerate}[i)]
\item The sum in \eqref{def:sigma} can be interpreted as a sum over all partitions of $n$ into $l$ distinct parts $\underline{u}_j$. The $v_j$ count the appearance of the parts $\underline{u}_j$. For example let $l=2$, $n=5$ and $r_1=r_2=1$ then we have five partitions of $5$ into $2$ distinct parts:
\begin{align*}
 5 &= \underline{4}+\underline{1} = \underline{3}+\underline{2} = \underline{3} + \underline{1} + \underline{1} 
 = \underline{2} +\underline{2} + \underline{1} = \underline{2} + \underline{1} + \underline{1} + \underline{1} \\
 &=  \underline{4}\cdot 1 +\underline{1}\cdot 1 = \underline{3} \cdot 1 + \underline{2} \cdot 1 = \underline{3} \cdot 1 + \underline{1} \cdot 2 = \underline{2}\cdot 2 + \underline{1} \cdot 1= \underline{2} \cdot 1 + \underline{1} \cdot 3 
\end{align*}
and therefore $\sigma_{0,0}(5) = 5$ and $\sigma_{2,1}(5) = 1^2\cdot 1^1 + 1^2\cdot 1^1 + 1^2\cdot 2^1 + 2^2 \cdot 1^1 + 1^2 \cdot 3^1 =11$\,. 

\item
The multiple divisor sum $\sigma_{\{0\}^l}$ counts the number of partitions of $n$ into $l$ distinct parts. Therefore the generating function of the partition functions $p(n)$ which counts all partitions of $n$ can be written as 
\[ \sum_{n>0} p(n) q^n = \sum_{l>0} [\{1\}^l] \,.\]
\item The brackets $[2,\dots,2]$ were already studied by P. A. MacMahon (see \cite{mm}) under the name of generalized divisor sums in the context of partitions. They were also studied in \cite{ar}where it was also shown, that they are quasi-modular forms. 
\end{enumerate}
\end{rem}

\begin{dfn}  Recall that for $s,z \in \C$, $|z|<1$ the polylogarithm $\Li_s(z)$ of weight $s$ is given by
\[ \Li_s(z) = \sum_{n>0} \frac{z^n}{n^s} \,. \]
We then define a \emph{normalized polylogarithm} by
\[ \Lit_{1-s}(z) := \frac{ \Li_{1-s}(z) }{\Gamma(s)}. \]
\end{dfn}

The normalized polylogarithm $\Li_{1-s}(z)$ extends to an entire function in $s$ and to a holomorphic function in $z$ where $|z|<1$. However for our purposes it is enough to know that for natural $s>0$ this is a rational function in $z$ with a pole at $z=1$ (c.f. \remref{rem_eulerpoly}). Now we can define brackets
as functions in $q$.

\begin{prop}
For $q \in \C$ with   $|q|<1$ and for all $s_1,\dots,s_l \in \N$
we can write the brackets as
\[ [s_1,\dots ,s_l] = \sum_{n_1 > \dots > n_l>0} \Lit_{1-s_1}\left(q^{n_1}\right) \dots \Lit_{1-s_l}\left(q^{n_l}\right) \,. \]
\end{prop}

\begin{prf} This follows directly from the definitions, see also \lemref{lem:eulerpol}.
\end{prf}

\begin{rem} 
\label{rem_eulerpoly}
As mentioned above the polylogarithms $\Li_{-s}(z)$ for $s \in \N$ are rational functions in $z$ with a pole in $z=1$. More precisely for $|z|<1$ they can be written as 
\[ \Li_{-s}(z) = \sum_{n>0} n^s z^n =\frac{z P_{s}(z)}{(1-z)^{s+1}}  \]
where $P_s(z)$ is the $s$-th Eulerian polynomial. Such a polynomial is given by
\[ P_s(X) = \sum_{n=0}^{s-1} A_{s,n} X^n \,, \]
where the  Eulerian numbers $A_{s,n}$ are defined by
\[ A_{s,n} = \sum_{i=0}^n (-1)^i \binom{s+1}{i} (n+1-i)^s\,.\]
Therefore the coefficients (the Eulerian numbers) of $P_s$ are positive. It fulfills the relation 
\[ P_{k+1}(t) = P_k(t) (1+kt) + t(1-t) P'_k(t)\, \]
and therefore $P_k(1)=k!$. 
For proofs of all these properties see for example \cite{df}. In particular the recursive formula can be found in \cite{df} as equation (3.3). \propref{prop:productli} then gives an expression for the product of Eulerian polynomials as rational linear combinations of polynomials in the form $(1-z)^j P_i(z)$ with $j,i \in \N$.  
\end{rem}

\begin{lem} \label{lem:eulerpol}
For $s_1,\dots,s_l \in \N$ we have
\[ [s_1,\dots,s_l] =  \frac{1}{(s_1-1)! \dots (s_l-1)!} \sum_{n_1 > \dots > n_l > 0} \prod_{j=1}^l \frac{q^{n_j} P_{s_j-1}\left( q^{n_j} \right)}{(1-q^{n_j})^{s_j}} \,, \]
where $P_k(t)$ is the $k$-th Eulerian polynomial. 
\end{lem}
\begin{prf} The claim follows directly  from \remref{rem_eulerpoly} because
\[\sum_{n_1 > \dots > n_l > 0} \prod_{j=1}^l \frac{q^{n_j} P_{s_j-1}\left( q^{n_j} \right)}{(1-q^{n_j})^{s_j}} = \sum_{n_1 > \dots > n_l > 0} \prod_{j=1}^l \sum_{v_j>0} v_j^{s_j-1} q^{v_j n_j} = \sum_{n>0} \sigma_{s_1-1,\dots,s_l-1}(n) q^n \,.   \] 
\end{prf}

The product of $[s_1]$ and $[s_2]$ can thus be written as
\begin{align*}
 [s_1] \cdot [s_2] &= \sum_{n_1 > n_2 > 0} \Lit_{1-s_1}\left(q^{n_1}\right) \Lit_{1-s_2}\left(q^{n_2}\right) + \sum_{n_2 > n_1 > 0} \dots + \sum_{n_1 = n_2 > 0} \Lit_{1-s_1}\left(q^{n_1}\right) \Lit_{1-s_2}\left(q^{n_1}\right) \\
 &= [s_1,s_2] + [s_2,s_1] + \sum_{n> 0} \Lit_{1-s_1}\left(q^n\right) \Lit_{1-s_2}\left(q^n\right) \,.
\end{align*}
In order to prove that this product is an element of $\filw_{s_1+s_2}(\MD)$ the product $\Lit_{1-s_1}\left(q^n\right) \Lit_{1-s_2}\left(q^n\right)$ must be  a rational linear combination of $\Lit_{1-j}\left(q^n\right)$ with $1 \leq j \leq s_1+s_2$. We therefore need the following

\begin{lem}  \label{lem_multpolylog}
For $a,b \in \N$ we have
\[ \Lit_{1-a}(z) \cdot  \Lit_{1-b}(z) = \sum_{j=1}^a \lambda^j_{a,b} \Lit_{1-j}(z) + \sum_{j=1}^b \lambda^j_{b,a}  \Lit_{1-j}(z)  + \Lit_{1-(a+b)}(z) \,, \]
where the coefficient $\lambda^j_{a,b}  \in \Q$ for $1 \leq j \leq a$ is given by
\[ \lambda^j_{a,b} = (-1)^{b-1} \binom{a+b-j-1}{a-j} \frac{B_{a+b-j}}{(a+b-j)!} \,. \]
\label{prop:productli}
\end{lem}
\begin{prf}
We prove this by using the generating function
\[ L(X) := \sum_{k>0} \Lit_{1-k}(z) X^{k-1} = \sum_{k>0} \sum_{n>0} \frac{n^{k-1} z^n}{(k-1)!} X^{k-1} = \sum_{n>0} e^{nX} z^n = \frac{e^X z}{1-e^X z} \,. \] 
With this one can see by direct calculation that
\[ L(X) \cdot L(Y) = \frac{1}{e^{X-Y}-1} L(X) + \frac{1}{e^{Y-X}-1} L(Y) \,. \]
By the definition of the Bernoulli numbers 
\[ \frac{X}{e^X-1} = \sum_{n\geq0} \frac{B_n}{n!} X^n  \]
this can be written as
\[ L(X) \cdot L(Y) = \sum_{n>0} \frac{B_n}{n!}(X-Y)^{n-1} L(X) + \sum_{n>0} \frac{B_n}{n!}(Y-X)^{n-1} L(Y) + \frac{L(X) - L(Y)}{X-Y}   \,. \]
The statement then follows by calculating the coefficient of $X^{a-1}Y^{b-1}$ in this equation.
\end{prf}

\begin{ex} We have $\lambda^1_{1,1} = B_1 = -\frac{1}{2}$ and thus 
\[ \Lit_{1-1}(z) \cdot  \Lit_{1-1}(z) = -\Lit_{1-1}(z) + \Lit_{1-2}(z) \,. \]
Therefore the product  $[1]\cdot[1]$ is given by
\[ [1]\cdot[1] = 2 [1,1] + [2] - [1] \,. \]
\end{ex}

More generally, \lemref{lem_multpolylog} implies the following explicit formula for the product in the length one case.

\begin{prop} \label{prop:l2-expli} We have the  formula
\begin{align*}
[s_1]\cdot [s_2]
=[s_1,s_2] +[s_2, s_1] + [s_1+s_2] 
+ 
  \sum_{j=1}^{s_1} \lambda^j_{s_1,s_2} [j] + \sum_{j=1}^{s_2} \lambda^j_{s_2,s_1} [j] \,.
  \end{align*} 
\begin{prf}This is a straightforward calculation \end{prf}
\end{prop}

In order to prove \thmref{prop:alg} we need to show that the above considerations work in general and not only in the length $1$ case. For this we use the notion of quasi-shuffle algebras (\cite{hi}). Let $A=\left\{ z_1 , z_2, \dots \right\}$ be the set of letters $z_j$ for each natural number $j\in \N$, $\Q A$ the $\Q$-vector space generated by these letters and $\Q\langle A \rangle$ the noncommutative polynomial algebra over $\Q$ generated by words with letters in $A$. For a commutative and associative product $\diamond$ on $\Q A$, $a,b \in A$ and $w,v \in \Q\langle A\rangle$ we define on $\Q\langle A\rangle$ recursively a product by $1*w=w*1=w$ and
\[ aw \ast bv := a(w \ast bv) + b(aw \ast v) + (a \diamond b)(w \ast v) \,. \]
Equipped with this product one has the
\begin{prop}\label{prop:HI-quasi}
The vector space $\Q\langle A\rangle$ with the product $\ast$ is a commutative $\Q$-algebra. 
\end{prop}
\begin{prf}
See \cite{hi} Theorem 2.1. 
\end{prf}

Motivated by the product expression of the polylogarithms in \lemref{prop:productli} we define the product $\diamond$ on $\Q A$ by
\[ z_a \diamond z_b =  \sum_{j=1}^a \lambda^j_{a,b} z_j + \sum_{j=1}^b \lambda^j_{b,a} z_j + z_{a+b}\,. \]
This is an commutative and associative product  on $\Q A$, because it arises from the product of the pairwise linearly independent polylogarithms $\Lit_{1-t}(z)$ in \propref{prop:productli},
and therefore   $\left( \Q\langle A\rangle, \ast \right)$ is a commutative $\Q$-algebra  
by \propref{prop:HI-quasi} above. 
\thmref{prop:alg} now follows from the next proposition.

\begin{prop}\label{prop:hoffmannalgebra}
For the linear map $[\, . \, ]: (\Q\langle A\rangle, *) \longrightarrow (\MD, \cdot)$ defined on the generators $w=z_{s_1}\dots z_{s_l}$ by $[w] := [s_1,\dots,s_l]$ we have
\[ [ w \ast v ]    = [w] \cdot [v]  \]
and therefore $\MD$ is a $\Q$-algebra and $[\, . \, ]$ an algebra homomorphism.
\end{prop}
\begin{prf}
This follows by the same argument as in the multiple zeta value case, see e.g. \cite{h1} Thm 3.2, by using induction on the length of the words $w$ and $v$ together with \propref{prop:productli}. 
\end{prf}

Now we have proven \thmref{thm:md-algebra}. As a special case of this theorem we have the following explicit formula.
\begin{ex} For $a,b,c \in \N$ we have
\begin{align*} 
[a]\cdot[b,c] &=
[z_a*z_bz_c] 
= [z_a z_b z_c + z_b z_a z_c +z_b z_c z_a + z_b(z_a \diamond z_c) + (z_a \diamond z_b)z_c ]\\
&=[a,b,c]+[b,a,c]+[b,c,a] + [a+b,c] + [b, a+c] \\
& \quad + \sum_{j=1 }^a \lambda_{a,c}^j[b,j]
 + \sum_{j=1 }^c  \lambda_{c,a}^j [b, j] 
 + \sum_{j=1 }^a \lambda_{a,b}^j [j, c]
+  \sum_{j= 1}^b \lambda_{b,a}^j [j, c].
\end{align*}
\end{ex}

We would like to point out another structure of the algebra $\MD$, which will be important later on when we consider the connection to multiple zeta values, and which was already mentioned in the introduction.

\begin{dfn}
We define the set of all admissible brackets $\MDA$ as the span of all brackets $[s_1,\dots,s_l]$ with $s_1 > 1$. With $\filwle_{k,l}(\MDA)$ we denote the admissible brackets of length $l$ and weight $k$ similar to the non-admissible case. 
\end{dfn}

With this we have the
 
\begin{thm}\label{thm:qmz-subalgebra}
The vector space $\MDA$ is a subalgebra of $\MD$.
\end{thm}
\begin{prf}
It is enough to show that $\MDA$ is closed under multiplication. Let $f=[a,\dots]$ and $g=[b,\dots]$ be elements in $\MDA$, i.e. $a>1$ and $b>1$. Due to \propref{prop:hoffmannalgebra} we have 
\[f \cdot g = [z_a w] \cdot [z_b v] = [ z_a w \ast z_b v ] \,,   \]
where $w,v \in \Q\langle A \rangle$ are words in the alphabet $A=\left\{ z_1 , z_2, \dots \right\}$. So in order to prove the statement we have to show that $z_a w \ast z_b v$ is a linear combination of words $z_c u \in \Q\langle A \rangle$ with $c>1$ and arbitrary words $u \in \Q\langle A \rangle$. By the definition of the quasi-shuffle product $\ast$ we have
\[ z_a w \ast z_b v = z_a(w \ast z_b v) + z_b(z_a w \ast v) + (z_a \diamond z_b)(w \ast v) \,. \]
The first two summands clearly fulfill this condition, because we assumed $a,b > 1$, so it remains to show that $z_a \diamond z_b \in \Q A$ is a linear combination of letters $z_j$ with $j>1$. Again by definition we obtain
\begin{align*}
 z_a \diamond z_b &=   z_{a+b} + \sum_{j=1}^a \lambda^j_{a,b} z_j + \sum_{j=1}^b \lambda^j_{b,a} z_j  \\
 &=  z_{a+b} + \left( \lambda^1_{a,b}+\lambda^1_{b,a} \right) z_1 + \sum_{j=2}^a \lambda^j_{a,b} z_j + \sum_{j=2}^b \lambda^j_{b,a} z_j \,,
\end{align*}
so it suffices to show that $\lambda^1_{a,b}+\lambda^1_{b,a}$ vanishes for $a,b > 1$. From the definition of $\lambda^j_{a,b}$ in \lemref{lem_multpolylog} it is easy to see that
\[ \lambda^1_{a,b}+\lambda^1_{b,a} = \left( (-1)^{a-1} + (-1)^{b-1} \right) \binom{a+b-2}{a-1} \frac{B_{a+b-1}}{(a+b-1)!} \,. \]
This term clearly vanishes when $a$ and $b$ have different parity. In the other case $a+b-1$ is odd and greater than $1$, as $a,b > 1$. It is well known that in this case $B_{a+b-1}=0$, from which we deduce that   $\lambda^1_{a,b}+\lambda^1_{b,a}=0$.
\end{prf}
 
\begin{thm} \label{thm:polyadalg}
\begin{enumerate}[i)]
\item We have  $\MD = \MDA[\,[1]\,]$. 
\item The algebra $\MD$ is a polynomial ring over $\MDA$ with indeterminate $[1]$, i.e. $\MD$ is isomorphic to $\MDA[\,T\,]$ by sending $[1]$ to $T$. 
\end{enumerate}
\end{thm}

\begin{prf} 
\begin{enumerate}[i)] 
\item First we show that any $f \in \filw_{k}(\MD)$ can be written as a polynomial in $[1]$.
If we show that for a fixed $l$ and $f \in \filwle_{k,l}(\MD)$ one can find 
$g_1 \in \filwle_{k,l}(\MDA)$
and $g_2,g_3 \in \filwle_{k,l-1}(\MD)$ such that $f$ can be written as
\begin{equation} \label{eq:polyrep}
f = g_1 + [1] \cdot g_2 + g_3 \,,
\end{equation}
then the claim follows directly by induction on $l$.

To show \eqref{eq:polyrep} it is clear that we can focus on the generators of $\MD$ which we write as $f=[ \{1\}^m , s_1, \dots , s_{l-m}]$, with $s_1 > 1$ and $k=m+s_1+\dots+s_{l-m}$. By induction over $m$ we prove that every element of such form can be written as in \eqref{eq:polyrep}. For $m=0$ it is $f \in \filw_{k}(\MDA)$, i.e. $g_1=f$ and $g_2=g_3=0$. For the induction step we obtain by the quasi-shuffle product
\begin{align*}
m \cdot [ \{ 1 \}^m , s_1, \dots, s_{l-m} ] &= [1] \cdot [\{1\}^{m-1},s_1,\dots,s_{l-m}] - g_3\\
&- \sum_{\substack{m_1+\dots+m_i = m \\ m_j \geq 0 \,, \forall j=1\dots i \\ m_1 < m}} [\{1\}^{m_1} , s_1, \{ 1 \}^{m_2} , \dots, s_{l-m} , \{ 1 \}^{m_{i}} ] \,.
\end{align*}
with $g_3 \in \filwle_{k,l-1}(\MD)$.
The elements in the sum start with at most $m-1$ ones, so we obtain a representation in the form of \eqref{eq:polyrep} inductively.

\item We have to show that $[1]$ is algebraically independent over $\MDA$ and therefore the
representation of $f\in \MD$ in i) as a polynomial in $[1]$ with coefficients in $\MDA$ is unique.  
From \propref{prop_Zkexpli}   we obtain that for $[s_1,\dots,s_l]\in \MDA$ with $s_1+\dots+s_l=k$ 
we have for $q$ close to $1$ the approximations $[s_1,\dots,s_l]\approx \frac{1}{(1-q)^k}$ and 
from \remref{rem:pupy} we know
$[1] \approx \frac{-\log(1-q)}{1-q}$. Therefore the only polynomial in $\MDA[T]$, which has $[1]$
as one of its roots,
is the constant polynomial $0$.
\end{enumerate}
\end{prf}

\begin{rem}
It is clear that $[1]$ is an irreducible element in the ring $\MD$, thus it is clear that
$\MD / \big([1]\cdot \MD\big)$ is a domain. But the non-obvious fact is that this domain can be represented by $\MDA$. 
\end{rem}
 
At the end of this section we want to mention two other subalgebras of $\MD$. For this denote by $\MDE$ the space spanned by  $1$ and all $[s_1,\dots,s_l]$ with $s_j$ even for all $0 \leq j \leq l$ and by $\MD^\sharp$ the space spanned by all by  $1$ and all $[s_1,\dots,s_l]$ with $s_j>1$.
\begin{prop}
$\MDE$ and $\MD^\sharp$ are subalgebras of $\MD$. 
\end{prop}
\begin{prf}
By the quasi-shuffle product formula \propref{prop:hoffmannalgebra} it is sufficient to show that for $1 \leq j \leq a$ the $\lambda^j_{a,b}  \in \Q$ given by
\[ \lambda^j_{a,b} = (-1)^{b-1} \binom{a+b-j-1}{a-j} \frac{B_{a+b-j}}{(a+b-j)!}  \]
vanish for $j$ odd if $a$ and $b$ are even to prove that $\MDE$ is a subalgebra of $\MD$. But this follows direclty by the fact that the $B_k$ vanish for odd $k>1$ and that the case $a+b-j=1$ does not occur since $j \leq a$ and $b\geq 2$. 

In order to prove that $\MD^\sharp$ is a subalgebra of $\MD$ we have to show that      
\[ \lambda^1_{a,b}+\lambda^1_{b,a} = \left( (-1)^{a-1} + (-1)^{b-1} \right) \binom{a+b-2}{a-1} \frac{B_{a+b-1}}{(a+b-1)!} \, \]
vanishes for $a,b > 1$. This term clearly vanishes when $a$ and $b$ have different parity. In the other case 
it is $B_{a+b-1}=0$, because
$a+b-1$ is odd and 
greater than $1$. Hence it is   $\lambda^1_{a,b}+\lambda^1_{b,a}=0$, whenever $a,b > 1$.
\end{prf}
 
The space $\MD^\sharp$ is studied further in \cite{BK}, where the authors consider a connection of this space to other $q$-analogues of multiple zeta values. 
 
\section{A derivation and linear relations in $\MD$}

Our strategy to prove \thmref{thm:derivative} is to use generating series of brackets. This allows
us to express the derivative in terms of elements in $\MD$. We make these calculations explicit in the case
of first in the length $1$ case and then  for the length $2$ case. Similar formulas for the general case are rather complicated.

\begin{lem}
The generating series $T(X_1,\dots,X_l)$ of brackets of length $l$ can be written as
\[  T(X_1,\dots,X_l) = \sum_{s_1,\dots,s_l >0} [s_1,\dots,s_l] X_1^{s_1-1} \dots X_l^{s_l-1}= \sum_{n_1,\dots,n_l>0} \prod_{j=1}^l \frac{e^{n_j X_j} q^{n_1+\dots+n_j}}{1-q^{n_1+\dots+n_j}} \,. \]

\end{lem}
\begin{prf}
This can be seen by direct computation using the geometric series and the Taylor expansion of the exponential function:
\begin{align*}
 \sum_{n_1,\dots,n_l>0} \prod_{j=1}^l \frac{e^{n_j X_j} q^{n_1+\dots+n_j}}{1-q^{n_1+\dots+n_j}}  =  &\sum_{n_1,\dots,n_l>0} \prod_{j=1}^l e^{n_j X_j} \sum_{v_j>0} q^{v_j(n_1+\dots+n_j)} \\
= &\sum_{n_1,\dots,n_l>0} \prod_{j=1}^l \sum_{k_j\geq0} \frac{n_j^{k_j}}{k_j!}  X_j^{k_j} \sum_{v_j>0} q^{v_j(n_1+\dots+n_j)} \\
\overset{u_j=v_j+\dots+v_l}{=} &\sum_{k_1,\dots,k_l\geq 0} \left( \sum_{\substack{u_1 > \dots>u_l>0 \\ n_1,\dots,n_l>0}} \frac{n_1^{k_1} \dots n_l^{k_l}}{k_1!\dots k_l!} q^{u_1 n_1 + \dots + u_l n_l} \right) X_1^{k_1} \dots X_l^{k_l} \\
=  &\sum_{s_1,\dots,s_l > 0} [s_1,\dots,s_l] X_1^{s_1-1} \dots X_l^{s_l-1}\,.
\end{align*}
\end{prf}

We now study  the derivative of brackets of length $1$, much of the formulas presented
for this purpose may implicitly found also in \cite{gkz}. In particular the next lemma is essentially a part of the calculation in the proof of Theorem 7 in \cite{gkz}. We give it nevertheless because it is a good preparation for the proof of our  \thmref{thm:derivative}. 

\begin{lem} \label{lem:shuffle2} \begin{enumerate}[i)]
\item
The product of two generating functions of multiple divisor sums of length $1$ is given by
\begin{align*}
T(X)\cdot T(Y) &= T(X+Y,X) + T(X+Y,Y) 
 					 - T(X+Y) + R_1(X,Y)\,
					 \intertext{where}
R_1(X,Y)&=\sum_{n>0} e^{n(X+Y)} \frac{q^n}{(1-q^n)^2}\,.
\end{align*}
 \item We have 
\[
\sum_{n>0} e^{nX} \frac{q^n}{(1-q^n)^2} = \sum_{k>0} \frac{\dif[k]}{k} X^{k}  +[2].
\]
In particular 
\[  R_1(X,Y) = \sum_{k>0} \frac{\dif[k]}{k} (X+Y)^{k} +[2] \,.\]
\end{enumerate}
\end{lem}
\begin{prf}  i) Remember that the generating functions are given by
\[ T(X) = \sum_{k>0} [k] X^{k-1} = \sum_{n>0} e^{nX} \qe{n} \] 
and 
\[ \qquad T(X,Y) = \sum_{s_1,s_2>0} [s_1,s_2] X^{s_1-1} Y^{s_2-1} = \sum_{n_1,n_2>0} e^{n_1X+n_2Y} \qe{n_1} \qe{n_1+n_2} \,.\]
With this in our hands we calculate
\begin{align*}
T(X) T(Y) &= \sum_{n_1,n_2 > 0} e^{n_1 X + n_2 Y} \qe{n_1} \qe{n_2} \\ 
&= \sum_{n_1>n_2>0}\dots+ \sum_{n_2 > n_1>0} \dots +  \sum_{n_2 = n_1>0} \dots =: F_1 + F_2 + F_3.
\end{align*}
For these terms we get furthermore
\begin{align*}
F_1 = &\sum_{n_1>n_2>0}  e^{n_1 X + n_2 Y} \qe{n_1} \qe{n_2} \\
 \overset{n_1 = n_2 + n_1'}{=} &\sum_{n_1',n_2>0}  e^{n_1' X + n_2 (X+Y)} \qe{n_1'+n_2} \qe{n_2} = T(X+Y,X)
\end{align*}
\begin{align*}
F_2 = &\sum_{n_2>n_1>0}  e^{n_1 X + n_2 Y} \qe{n_1} \qe{n_2} \\
 \overset{n_2 = n_1 + n_2'}{=} &\sum_{n_1,n_2'>0}  e^{n_1 (X+Y) + n_2 Y} \qe{n_1} \qe{n_1+ n_2'} = T(X+Y,Y).
\end{align*}

Using $\left( \qe{n} \right)^2 = \frac{q^n}{(1-q^n)^2} - \frac{q^n}{1-q^n}$, we get for the last term
\begin{align*}
F_3 = &\sum_{n_1=n_2 > 0}  e^{n_1 (X +Y)} \left( \qe{n_1} \right)^2 \\
= &\sum_{n > 0}  e^{n (X+Y)}  \frac{q^n}{(1-q^n)^2} - \sum_{n > 0}  e^{n (X+Y)}  \frac{q^n}{(1-q^n)} \\
= &\, R_1(X,Y) - T(X+Y) \,.
\end{align*}

ii) This can be seen by direct computation. First observe
\[ \dif T(X) = \sum_{k>0} \dif[k] X^{k-1}  = \dif \sum_{n>0} e^{nX} \qe{n} = \sum_{n>0} n e^{nX} \frac{q^n}{(1-q^n)^2} \,.\]
and then use this to evaluate
\begin{align*}
 \sum_{k>0} \frac{\dif[k]}{k} X^{k} &= \sum_{k>0} \int_0^X \dif[k] t^{k-1} dt \\
&=  \int_0^X  \dif T(t) dt = \sum_{n>0} \int_0^X n e^{nt} dt \frac{q^n}{(1-q^n)^2} \\
&= \sum_{n>0} e^{nX} \frac{q^n}{(1-q^n)^2} - \sum_{n>0} \frac{q^n}{(1-q^n)^2}= \sum_{n>0} e^{nX} \frac{q^n}{(1-q^n)^2} - [2] \,.
\end{align*}
\end{prf}

We now want to give explicit expressions for the derivative of multiple divisor sums of length $1$, which follow from the  lemmas  above:
\begin{prop} \label{cor:formularfordk}
For $s_1,s_2$ with $s_1+s_2>2$ and $s=s_1+s_2-2$ we have the following expression for $\dif[s]$:
\[ \binom{s}{s_1-1} \frac{\dif[s]}{s} = [s_1]\cdot [s_2] +\binom{s}{s_1-1} [s+1] - \sum_{a+b=s+2} \left( \binom{a-1}{s_1-1}+\binom{a-1}{s_2-1} \right) [a,b] \,.\]
\end{prop}
\begin{prf}
This is a direct consequence of \lemref{lem:shuffle2} by considering the coefficient of $X^{s_1-1} Y^{s_2-1}$ in the equation 
\begin{align*}
 T(X)\cdot T(Y) &= T(X+Y,X) + T(X+Y,Y)- T(X+Y) 
 					+  \sum_{k>0} \frac{\dif[k]}{k} (X+Y)^{k} +[2]  \,.
\end{align*}
by using 
\begin{align*}
T(X+Y,X) + T(X+Y,Y) &=\sum_{\substack{s_1,s_2 >0\\ a+b=s_1+s_2}} \left(\binom{a-1}{s_1-1}+\binom{a-1}{s_2-1} \right) [a,b] X^{s_1-1} Y^{s_2-1} \,,\\
T(X+Y) &= \sum_{s_1,s_2 >0} \left( \binom{s_1+s_2-2}{s_1-1} [s_1+s_2-1] \right) X^{s_1-1} Y^{s_2-1} \,, \\
\sum_{k>0} \frac{\dif[k]}{k} (X+Y)^{k} &= \sum_{s_1,s_2 >0} \left(  \binom{s_1+s_2-2}{s_1-1} \frac{d[s_1+s_2-2]}{s_1+s_2-2}\right) X^{s_1-1} Y^{s_2-1} \,.
\end{align*}
\end{prf}

\begin{ex} In the following formulas we used the explicit description for the product given in 
\propref{prop:l2-expli}.
\begin{enumerate}[i)]
\item In the smallest case $s=1$ there is just one choice given by $s_1=1$, $s_2=2$:  
\[ \dif[1] = [3] + \frac{1}{2}[2]-[2,1] \,. \]
\item For $s=2$ we can choose $s_1=1, s_2 =3$ and $s_1=s_2=2$ and therefore we get the two expressions:
\begin{align*}
\dif[2] &= 2[4] + [3] +\frac{1}{6} [2] -2[2,2] - 2[3,1] \,,\\
\dif[2] &= [4] +2 [3] - \frac{1}{6} [2] - 4[3,1]  \,,
\end{align*}
from which the first linear relation in weight $4$ follows:
 \[ [4] =  2[2,2] - 2 [3,1] + [3] - \frac{1}{3} [2] \,.\]
\item In the case $s=3$ one again gets two expressions and therefore one relation. 
\item For $s=4$ one has $s_1=1$, $s_2=5$ or $s_1=2$, $s_2=4$ and $s_1=s_2=3$ which gives
\begin{align*}
\dif[4] &= 4 [6] +2 [5] + \frac{1}{3} [4] - \frac{1}{180} [2] -4 [2,4] -4 [3,3]-4 [4,2]-4 [5,1] \,, \\
\dif[4] &= [6] +4 [5]-\frac{1}{12} [4]+ \frac{1}{180} [2] -2 [3,3]-3 [4,2]-8 [5,1] \,, \\
\dif[4] &= \frac{2}{3} [6] +4 [5] - \frac{1}{180} [2] -4 [4,2]-8 [5,1] \,.
\end{align*}
From which the following two relations follow
\begin{align*}
5 [6] &= 3 [5]-\frac{1}{2} [4] +6 [2,4]+6 [3,3]-6 [5,1] \\
3 [6] &= 2 [5]-\frac{5}{12}  [4] + \frac{1}{90} [2] +4 [2,4]+2 [3,3]+[4,2]-4 [5,1] \,.
\end{align*}
\end{enumerate}
\end{ex}

\begin{thm} \label{thm:dim2} Suppose $k \ge 4$, then there are at least $\lfloor \frac{k}{2} \rfloor-1$ linear relations in the
generators of $\grwl_{k,2} ( \MDA)$.  
\end{thm}
\begin{prf}
It is clear that the expressions for $\dif[k-2]$ in \propref{cor:formularfordk} are symmetric in $s_1$ and $s_2$. There are $\lfloor \frac{k}{2} \rfloor$ choices for $s_1$ and $s_2$ with $s_1+s_2 = k$ and 
$s_1 \leq s_2$.  For each such choice we get a 
  different expression for $\dif[k-2]$, because for $s_1 \leq s_2$ it only contains the length $2$ terms $[s_1+1,s_2-1],\dots, [s_1+s_2-1,1] \in \MDA$ with non vanishings coefficients. This can be seen if we rewrite the statement by using the stuffle product $[s_1] \cdot [s_2] = [s_1,s_2] + [s_2,s_1] + [s_1 \diamond s_2]$:  

\begin{align*}
&\binom{k-2}{s_1-1} \frac{\dif[k-2]}{k-2} \\
&= [s_1 \diamond s_2]  +\binom{k-2}{s_1-1} [k-1] - \sum_{\substack{a+b=k \\ a > s_1}} \left( \binom{a-1}{s_1-1}+\binom{a-1}{s_2-1} - \delta_{a,s_2} \right) [a,b]  \,. 
\end{align*}
By the same considerations as in proof of \thmref{thm:qmz-subalgebra} we find that $[s_1 \diamond s_2] \in \MDA$. Therefore we get $\lfloor \frac{k}{2} \rfloor-1 $ relations. 
\end{prf}

We have checked that for $k\le 20$ we get all relations in length two by the above method, cf. \thmref{thm-dklmda}. This give some evidence for
\begin{conj}
For all weights $k \ge 4$ the number of linear relations in the generators of 
$\grwl_{k,2}\MDA$ equals  $\lfloor \frac{k}{2} \rfloor-1 $.    
\end{conj}

Now we want to consider the derivative in the length two case.

\begin{lem} \label{lem:len2v2}
 The product of two generating functions of multiple divisor sums of length $1$ and $2$ is given by
\begin{align*}
 T(X)\cdot T(Y,Z) =& T(X+Y,Y,Z) + T(X+Y,X+Z,Z) + T(X+Y,X+Z,X) \\
 					&  -T(X+Y,Z) - T(X+Y,X+Z) + R_2(X,Y,Z)  \,,
\end{align*}
where 
\begin{align*}
R_2(X,Y,Z) =&  
\sum_{n_1,n_2>0} e^{n_1(X+Y)+n_2Z} \frac{q^{n_1}}{(1-q^{n_1})^2} \qe{n_1+n_2} \\ 
 					&+\sum_{n_1,n_2>0} e^{n_1(X+Y)+n_2(X+Z)} \qe{n_1} \frac{q^{n_1+n_2}}{(1-q^{n_1+n_2})^2}
\end{align*}
\end{lem}

\begin{prf} i)
We again split the sum into different parts as in the case for $T(X)T(Y)$:
\begin{align*}
T(X)\cdot T(Y,Z) &= \sum_{n_1,n_2,n_3 > 0} e^{n_1X+n_2Y+n_3Z} \qe{n_1}\qe{n_2}\qe{n_2+n_3} \\
&= \sum_{n_2 > n_1} \dots + \sum_{n_1 > n_2 +n_3} \dots + \sum_{n_2+n_3 > n_1 > n_2} \dots + \sum_{n_1=n_2} \dots + \sum_{n_1=n_2+n_3}\dots \\
&=: F_1 + F_2 + F_3 + F_4 + F_5 \,.
\end{align*}
The proof of $F_1 + F_2 + F_3 = T(X+Y,Y,Z) + T(X+Y,X+Y,Z) + T(X+Y,X+Y,X)$ is similar to the calculation in the lemma above and we leave it out here.
The evaluation of $F_4$ and $F_5$ are similar and we therefore just illustrate the $F_4$ case:
\begin{align*} 
F_4 &= \sum_{n_1=n_2, n_3 > 0} e^{n_1(X+Y)+n_3Z} \left(\qe{n_1}\right)^2 \qe{n_1+n_3} \,, \\
&= \sum_{n_1=n_2, n_3 > 0} e^{n_1(X+Y)+n_3 Z} \left(  \frac{q^{n_1}}{(1-q^n)^2} - \qe{n_1} \right) \qe{n_1+n_3} \,, \\
&= \sum_{n_1, n_3 > 0} e^{n_1(X+Y)+n_3 Z}  \frac{q^{n_1}}{(1-q^n)^2} \qe{n_1+n_3} - T(X+Y,Z) \,.
\end{align*}

\end{prf}

\begin{dfn}
We define  the operator $D(f)$ on functions in $X$ by  
\[ D(f) = \left( \frac{\partial}{\partial X} f \right) \Big|_{X=0}. \]
\end{dfn}

Observe that $D(R_1(X,Y)) = \dif T(Y)$ and for the length $2$ it holds

\begin{lem} \label{lem:len2v2-r2}
We have 
\begin{align*}
D( R_2(X,Y,Z))  
= \dif T( Y, Z). 
\end{align*}
\end{lem}

 \begin{prf}
For the two summands of $R_2(X,Y,Z)$ one gets
\begin{small}
\begin{align*}
D\left( \sum_{n_1,n_2>0} e^{n_1(X+Y)+n_2Z} \frac{q^{n_1}}{(1-q^{n_1})^2} \qe{n_1+n_2} \right) &=\\   \sum_{n_1,n_2>0} n_1 e^{n_1 Y+n_2Z}& \frac{q^{n_1}}{(1-q^{n_1})^2} \qe{n_1+n_2}\,, \\
D\left( \sum_{n_1,n_2>0} e^{n_1(X+Y)+n_2(X+Z)} \qe{n_1} \frac{q^{n_1+n_2}}{(1-q^{n_1+n_2})^2} \right) &=\\ \sum_{n_1,n_2>0} (n_1+n_2) e^{n_1 Y+n_2 Z}& \qe{n_1} \frac{q^{n_1+n_2}}{(1-q^{n_1+n_2})^2}\,.
\end{align*}
\end{small}\\
Adding these two terms one obtains $\dif T(Y,Z)$, because with $\dif \qe{n} = \frac{ n \cdot q^n}{(1-q^n)^2}$ and the product formula we obtain
\begin{align*}
&\dif T(Y,Z) = \dif \sum_{n_1,n_2>0} e^{n_1Y+n_2Z} \qe{n_1} \qe{n_1+n_2} =\\
&\sum_{n_1,n_2>0} \left(  n_1 e^{n_1 Y+n_2Z} \frac{q^{n_1}}{(1-q^{n_1})^2} \qe{n_1+n_2} + (n_1+n_2) e^{n_1 Y+n_2 Z} \qe{n_1} \frac{q^{n_1+n_2}}{(1-q^{n_1+n_2})^2} \right) 
\end{align*}
\end{prf}

\begin{prop} \label{prop:len2}
The derivative of $[s_1,s_2]$ can be written as
\begin{align*}
\dif&[s_1,s_2] = [2]\cdot [s_1,s_2] - s_1 [s_1+1,s_2,1] - s_2 [s_1, s_2+1,1] - [s_1,s_2,2]\\
&- \left( \sum_{a+b=s_1+2} (a-1) [a,b,s_2] + \sum_{a+b=s_2+1} s_1 [s_1+1,a,b] + \sum_{a+b=s_2+2} (a-1) [s_1,a,b] \right)\\
&+ 2 s_1 [s_1+1,s_2] + s_2 [s_1,s_2+1] \,.\\
\end{align*}
\end{prop}
\begin{prf}
This follows directly from \lemref{lem:len2v2} by applying the operator 
\[ D(f) = \left( \frac{d}{dX} f \right) \big|_{X=0} \]
on both sides of the equation. 
It is straightforward to calculate $D(T(\dots,\dots))$ for the various generating series $T(\dots,\dots)$
in  \lemref{lem:len2v2}, e.g.,
the lefthand side becomes $[2] \cdot T(Y,Z)$.  
By means of \lemref{lem:len2v2-r2} the claim follows easily by  collecting all the terms. 
\end{prf}

\begin{ex}
\begin{enumerate}[i)]
\item For $s_1=s_2=1$  
 \propref{prop:len2} gives the  representation of $\dif[1,1]$ already mentioned above in 
\eqref{example:derivative4}: 
\[ \dif [1,1] =[3,1] +[1,3]+ \frac{3}{2}[2,1] + \frac{1}{2}[1,2]- 2 [2,1,1]-[1,2,1] \,; \]
here we used the quasi-shuffle product 
\[ [2]\cdot [1,1] =  [3,1] + [1,3] + [2,1,1]+[1,2,1]+[1,1,2] - \frac{1}{2} [2,1] - \frac{1}{2} [1,2] \,.\]

\item  For $s_1=1$, $s_2=2$ the corollary gives
\[ \dif[1,2] = -\frac{1}{6} [1,2] +2 [1,3]+[1,4]+\frac{3}{2} [2,2] +[3,2]-4 [1,3,1]-[2,1,2]-2 [2,2,1] \,. \] 

\item  For $s_1=1$, $s_2=2$ the corollary gives
\[ \dif[2,1] =  -\frac{1}{6} [2,1] + \frac{1}{2} [2,2] +[2,3]+4 [3,1]+[4,1]-[2,2,1]-6 [3,1,1] \,.\]

\item The case $s_1=s_2=2$ is given by
\[ \dif[2,2] = -\frac{1}{3} [2,2] +2 [2,3]+[2,4]+4 [3,2]+[4,2]-4 [2,3,1]-4 [3,1,2]-4 [3,2,1] \,. \]
\end{enumerate}

\end{ex}

At this point we like to indicate that the Leibniz rule is another source of linear relations in $\MD$.

\begin{ex}
\label{ex:rel5+6}
\begin{enumerate}[i)]
\item By means of the Leibniz rule and the quasi-shuffle product we have
\[  \dif[1] \cdot [2] + [1] \cdot \dif[2]
  = \dif( [1] \cdot [2]) =  \dif\left( [1,2] + [2,1] + [3] - \frac{1}{2} [2] \right) \,. \]
Evaluating both sides separately we deduce the following linear relation in length $3$
\begin{align}\label{eq:rel-w5}
 [5] =&  2 [3,1,1]-[2,2,1] +[2,3] + 2 [3,2]-[4,1] \notag \\
 			&+ \frac{1}{2} [4] + \frac{1}{2} [2,2] - 2 [3,1] + \frac{1}{6} [2,1]-\frac{1}{12} [2] + \frac{1}{12} [3] \,.
\end{align}
\item Using the same argument for $[1] \cdot [3]$ we have
\[ \dif[1] \cdot [3] + [1] \cdot \dif[3] = \dif( [1] [3] ) = \dif\left( [1,3] + [3,1] + [4] + \frac{1}{12} [2] - \frac{1}{2} [3] \right)  \]  
from which the following relation in weight $6$ follows
\begin{align*}
[6] = &\frac{1}{120} [2] - \frac{1}{24} [3]+  \frac{1}{2} [5]+\frac{1}{4}  [2,2]-[2,2,2]+ \frac{1}{2} [2,3]-[2,3,1]\\
&+[2,4]+ \frac{1}{12} [3,1]+2 [3,1,2]-[3,2]-3 [4,1]+3 [4,1,1]+5 [4,2]-[5,1] \,.
\end{align*}
\end{enumerate}
\end{ex}

\begin{thm}\label{thm:dim3} 
\begin{enumerate}[i)]
\item There is a linear relation in the generators of $\grwl_{5,3}(\MDA)$.
\item There are at least $3$ linear relations in the generators of $\grwl_{6,3}(\MDA)$.
 \end{enumerate}
\end{thm}
\begin{prf} i) From \exref{ex:rel5+6} i) we deduce the relation 
$$ 0 \equiv 2 [3,1,1] - [2,2,1]$$ 
in $\grwl_{5,3}(\MDA)$.

ii) From \exref{ex:rel5+6} ii) we deduce the relation 
$$0 \equiv 2 [3,1,2]+3 [4,1,1]$$ 
in $\grwl_{6,3}(\MDA)$ and from \exref{ex:oldrelations} we deduce  
that 
\begin{align*}
0 &\equiv 3 [2,2,2] - [3,1,2]- [2,3,1]-[3,2,1]\\
0& \equiv -4 [2,3,1]-4 [3,1,2]- [3,2,1]+9 [4,1,1]
\end{align*}
\end{prf}
 
We finally want to prove that the map $\dif$ is a derivation for arbitrary length using the same combinatorial arguments as in the length one and two cases but without calculating explicit representations for $\dif[s_1,\dots,s_l]$. 

\begin{prf}(of \thmref{thm:derivative})
To prove this statement we are going to use the same combinatorial arguments as in the \lemref{lem:len2v2}, \lemref{lem:len2v2-r2} and \propref{prop:len2} in a general way which means that we have 
 \begin{align} \label{eq:deriveq}
 \begin{split}
T(X) \cdot &T(Y_1,\dots,Y_l) =
\sum_{m,n_1,\dots,n_l>0} e^{mX + n_1Y_1 + \dots + n_l Y_l} \qe{m} \qe{n_1} \dots \qe{n_1+\dots+n_l} \\
=&T(X+Y_1,\dots,X+Y_l,X) + \sum_{j=1}^l T(X+Y_1,\dots,X+Y_j,Y_j,\dots,Y_l)   \\
&+R_l - \sum_{j=1}^l T(X+Y_1,\dots,X+Y_j,Y_{j+1},\dots,Y_l)\,,
\end{split}
\end{align}
where 
\[ R_l = \sum_{j=1}^l \left( \sum_{n_1,\dots,n_l >0} e^{n_1(X+Y_1) + \dots + n_j (X+Y_j) + n_{j+1} Y_{j+1} + \dots + n_l Y_l} \prod_{i=1}^l \frac{q^{n_1+\dots+n_i}}{(1+q^{n_1+\dots+n_i})^{\delta_{i,j}+1}}\right) \,.\]
This can be seen by splitting up the sum in the same way as above. The first line comes from the parts where one sums over the ordered pairs $n_1 + \dots + n_{j-1} < m < n_1 + \dots + n_j$ for $j=1,\dots ,l$ and $n_1 + \dots + n_l < m$. Setting $m = n_1 + \dots + n_{j-1} + m'$ and $n_j = m' + n_j'$ for these terms it is easy to see that one gets the sum over $m',n_1,\dots,n_j',\dots,n_l$ which then gives $T(X+Y_1,\dots,X+Y_j,Y_j,\dots,Y_l)$.\\
The second line arises from the sum over $m = n_1 + \dots + n_j$. In this case one again uses the identity  
\[ \left( \qe{n} \right)^2 = \frac{q^n}{(1-q^n)^2} - \frac{q^n}{1-q^n}\]
 from which the rest follows easily. \\
Letting the operator $D(f) = \left( \frac{d}{dX} f \right) \big|_{X=0}$ act on this it is easy to see that the last term then becomes
\[ D(R_l) = \sum_{j=1}^l \left( \sum_{n_1,\dots,n_l >0} (n_1+\dots+n_j) e^{n_1 Y_1 + \dots + n_l Y_l} \prod_{i=1}^l \frac{q^{n_1+\dots+n_i}}{(1+q^{n_1+\dots+n_i})^{\delta_{i,j}+1}}\right) \] 
and this is exactly $\dif T(Y_1,\dots,Y_l)$ which can be seen by induction on $l$ and the product formula.
The product on the left becomes $[2] T(Y_1 , \dots , Y_l)$ and the remaining terms on the right all have elements in $\filwle_{k+2,l+1}(\MD)$ as their coefficients and therefore the statement follows. 
\end{prf}

\begin{prop} \label{prop:qmzcloseddif}
The space $\MDA$ is closed under $\dif$.
\end{prop}
\begin{prf}
This follows directly by the proof of \thmref{thm:derivative} since in the formula for $\dif T(Y_1,\dots,Y_l)$, which one obtains by applying $D$ to equation \eqref{eq:deriveq}, it is easy to see that the coefficients of the monomials which contains a $Y_1$ are all in $\MDA$. 
\end{prf}
 
\begin{rem}
We didn't give an explicit formula for the derivative of brackets of length $l$, since a general formula seems to be confusing. But for a specific bracket one can get its derivative  by applying  first the operator $D$ to the equation \eqref{eq:deriveq} and then collecting the corresponding coefficients. For example for $l=3$ one can deduce
\begin{align*}
\dif [2,1,1] &= -\frac{1}{6} [2,1,1] + \frac{1}{2} [2,1,2] -[2,1,2,1]+[2,1,3]+ \frac{3}{2} [2,2,1] \\
 						 &-2\, [2,2,1,1]+[2,3,1]+6 [3,1,1]-8 [3,1,1,1]+[4,1,1]. 
\end{align*} 
\end{rem}

\begin{rem} Changing the perspective we can view \thmref{thm:derivative})
and its special cases  \lemref{lem:len2v2}, \lemref{lem:len2v2-r2}  and \propref{prop:len2}
as  results, which express the failure of the 
shuffle relation for $[s]\cdot[s_1, \ldots, s_l]$ in terms of multiple
 divisor functions of lower weight and length and derivatives. An optimistic guess
 is that this is also the case for more complicated products. We want to come back to this in \cite{BBK}.
 \end{rem}
 
\section{The subalgebra of (quasi-)modular forms}

We call 
\[ G_{k} =  \frac{\zeta(k)}{(2 \pi i)^k} + \frac{1}{(k-1)!} \sum_{n>0} \sigma_{k-1}(n) q^n  =     \frac{\zeta(k)}{(2 \pi i)^k}   + [k] \,. \]
the Eisenstein series of weight $k$. For even $k=2n$ due to Eulers theorem we have in addition
\[ \zeta(2n) = \frac{(-1)^{n-1} B_{2n} (2\pi)^{2n}}{2(2n)!}  \]
and therefore 
\[ G_{2n} =   -\frac{1}{2} \frac{B_{2n}}{(2n)!}  + [2n] 
\in \filw_{2n}(\MD),\]
for example
\[ G_2 = -\frac{1}{24} + [2] \,,\quad G_4 = \frac{1}{1440} + [4] \,, \quad G_6 = -\frac{1}{60480} + [6] \,.\]

\begin{prop} \label{prop:mf-neu} 
\begin{enumerate}[i)]
\item The ring of modular forms $M(\Gamma_1)$ for $\Gamma_1 = SL_2(\Z)$ 
and the ring of quasi-modular forms 
$\widetilde{M}(\Gamma_1)$ are graded subalgebras of $\MD$.
\item The $\Q$-algebra of quasi-modular forms   
$\widetilde M_k(\Gamma_1)$ 
is closed under the derivation $\dif$ and therefore it is a subalgebra of the graded differential algebra 
$(\MD, \dif)$. 

\item We have the following inclusions of $\Q$-algebras
\[ M_k(\Gamma_1) \subset  \widetilde{M}(\Gamma_1) \subset \MDE \subset \MD^\sharp \subset \MDA \subset \MD \,. \]  
\end{enumerate}
\end{prop}

\begin{prf} Let $M_k(\Gamma_1)$ (resp.  $\widetilde M_k(\Gamma_1)$) be the space of (quasi-)modular forms of weight $k$ for $\Gamma_1$. Then the first claim follows directly from the well-known facts
\begin{align*}
M(\Gamma_1) &= \bigoplus_{k>1} M(\Gamma_1)_k = \Q[G_4,G_6]\\
\widetilde{M}(\Gamma_1) &=\bigoplus_{k>1} \widetilde{M}(\Gamma_1)_k = \Q[G_2,G_4,G_6]\,.
\end{align*} 
 The second claim is a well known fact in the theory of quasi-modular forms and a proof can be found in \cite{dz} p. 49. It suffices to show that the
 derivatives of the generators are given by
\begin{align*} 
\dif G_2 &=\dif [2] = 5 G_4 - 2 G_2^2 \,,\qquad \dif G_4 = 15 G_6 - 8 G_2 G_4 \,,  \\
\dif G_6 &= 20 G_8 - 12 G_2 G_6 = \frac{120}{7} G_4^2 - 12 G_2 G_6 \,. 
\end{align*} 
The last statement follows immediately by i) and the results before. 
\end{prf}

\begin{rem}   
The above formulas for $\dif [2], \dif [4]$ and $\dif [6]$ can also be proven with \propref{cor:formularfordk}.
\end{rem}
\begin{ex} The theory of modular forms yield linear relations in $\MD$. We indicate
here how to derive such a relation in weight $8$.  
It is a well-known fact from the theory of modular forms that $G_4^2 = \frac{7}{6} G_8$ because the space of weight $8$ modular forms is one dimensional. We therefore have   
\begin{align*}
  \frac{1}{2073600}  + \frac{1}{720} [4]  + [4] \cdot [4] 
&=
\left( \frac{1}{1440} 
+ [4] \right)^2 =     \frac{1}{2073600} + \frac{7}{6} [8]\,.
\end{align*}
Using the quasi-shuffle product from \propref{prop:productli} we get
\begin{align*}
 [4] \cdot [4] 
&= 2 [4,4] + [8] + \frac{1}{360} [4]- \frac{1}{1512} [2] \,,
\end{align*}
which then gives the following relation in weight $8$:
\[ [8] = \frac{1}{40} [4] - \frac{1}{252} [2] + 12 [4,4]\,.  \]
\end{ex}

It is well known that the weight
is additive for multiplication of modular forms.
The above  relation shows that the length is
not additive with respect to the
multiplication of modular forms. 
 
\begin{prop} The algebra of modular forms is graded with respect to the weight
and filtered with respect to the length. We have
\begin{align*}
\sum_k \dim_\Q \grwl_{k,l} M(\Gamma_1)\, x^k y^l &=
 1 + \frac{x^4}{1-x^2} y + \frac{x^{12}}{(1-x^4)(1-x^6)}y^2\,,\\
\intertext{in particular} 
\sum_k \dim_\Q M_k(\Gamma_1)\, x^k &= \frac{1}{(1-x^4)(1-x^6)}\,.\\
\end{align*}
\end{prop}

\begin{prf}   For each $k$ there is an Eisenstein series $G_k$ and this is the only element
of length $1$ in $M_k(\Gamma_1)$. Now
the first statement follows immediately from the fact that 
the polynomials $G_a G_b$ with $a+b=k$ generate $M_k(\Gamma_1)$ as an vector space \cite{zag77}. 
Setting $y$ in the first formula we see again that
the modular forms $G_4$ and $G_6$ generate $M(\Gamma_1)$ as an algebra. 
\end{prf}

Notice that because of  \thmref{thm:lowdim_fil_md} we know all relations in $\filwle_{8,2}(\MD)$ 
and therefore we could give a purely algebraic proof the relation $G_4^2 = \frac{7}{6} G_8$ 
without using the theory of modular forms, which relies on complex analysis.  
Moreover, again using  \thmref{thm:lowdim_fil_md}, we can prove in $\filw_{12}(\MD)$ 
new identities for the cusp form $\Delta = \sum_{n>0} \tau(n) q^n$. 

\begin{prop}\label{prop:deltarep} For $(a,b) \in \{(2,4),(4,6),(6,8),(8,10),(10,11),(11,12)\}$ the cusp form $\Delta \in  S_{12}$ can be uniquely written as
\[ \Delta = \frac{2^b + 50}{2^{b}-2^{a}} \cdot [a] + \frac{2^a + 50}{2^{a} - 2^{b}}\cdot [b] + \sum_{m+n=12} d_{m,n} \cdot [m,n] \,, \]
where $d_{m,n} \in \Q$. Moreover, any other representation of $\Delta$ in $\filwle_{12,2}(\MD)$
is a linear combination of these six representations.
\end{prop}
\begin{prf} By  \thmref{thm:lowdim_fil_md} we just have to solve systems of linear equations coming from the coefficients of the brackets in question. Using
 the relations coming from \propref{cor:formularfordk} this can be made very efficient with the
 computer.
\end{prf}

Taking a suitable linear combination of the identities in \propref{prop:deltarep} we get 
the representation \eqref{eq:deltal2} of $\Delta$ given in the introduction.


\begin{rem}
At the end of this section we just want to give a short remark concerning the arithmetical aspect of the relations in \propref{prop:deltarep} on which we don't want to focus in detail in these notes.
Formulas like the ones above give several representation of the Fourier coefficients of cusp forms in terms of multiple divisor sums. One can also see the well-known congruence $\tau(n) \equiv \sigma_{11}(n) \mod 691$ and it is easy to derive a lot of other congruences involving $\tau(n)$ and the brackets out of such relations. 
\end{rem}

\section{Experiments and conjectures: dimensions}

%
%
%
%
In this section we present data of some computer calculations regarding the number of linear independent brackets with length and weight smaller or equal to $15$.  In some cases we can prove these
bounds to be sharp. Based on these experiments, we make a conjecture on the dimension
of the graded pieces of $\MDA$ and therefore also for $\MD$.
We first recall our results on the algebraic structure of $\MD$ and $\MDA$, where $\MDA$ is the sub algebra of $\MD$ generated by admissible brackets.  
Both are a bi-filtered algebras with respect to the
filtration $\filw_\bullet$ given by the weight and the filtration $\fille_\bullet$ given by the
length. 
Therefore as   vector spaces we have
\begin{align}
\MD &\cong \bigoplus_{k} \grw_{k}(\MD) \cong \bigoplus_{k}\bigoplus_{l\le k} \grwl_{k,l}(\MD) \label{eq:dirsum_md}\\
\MDA &\cong \bigoplus_{k} \grw_{k}(\MDA) \cong \bigoplus_{k} \bigoplus_{l\le k-1}\grwl_{k,l}(\MDA).
\label{eq:dirsum_mda}
\end{align}
  
\begin{prop} \label{prop:grad-fin} In the direct sums in \eqref{eq:dirsum_md} and \eqref{eq:dirsum_mda} each summand is a finite dimensional vector space. In particular, we have  
\begin{align*} 
\dim_\Q \grwl_{k,l}( \MD )  \le \binom{k-1}{l-1}
\,,\qquad 
\dim_\Q \grwl_{k,l} ( \MDA )  \le \binom{k-2}{l-1}.
\end{align*}
\end{prop}
\begin{prf}
Let $b(k,l)$ denote the number of brackets $[s_1,\dots,s_l]$ of weight $k$ and length $l$, i.e. $s_1+\dots+s_l = k$ and let $a(k,l)$ denote the number of admissible brackets of this type, i.e. $s_1+\dots+s_l = k$ with $s_1>1$. 
It suffices to show
\begin{align}\label{eq:def_bkl_akl}
 b(k,l) = \binom{k-1}{l-1} \,,\qquad a(k,l) = \binom{k-2}{l-1} \,. 
\end{align}
Now, if we write $k=1+\dots+1$, then these formulas
are an easy combinatorial fact, which can be seen by counting the possible ways of replacing $l-1$ of $k-1$ plus symbols by a semi-column    and then interpreting the remaining sums as tuples $(s_1,\dots,s_l)$ 
(resp. $k-2$ since we can't replace the first plus symbol). 
\end{prf}

\begin{dfn}
We define 
\begin{align*}
d'(k,l)&= \dim_\Q  \grwl_{k,l}( \MDA ). 
\end{align*}
\end{dfn}

The next proposition shows that, in order to understand the dimensions of the  various subspaces of $\MDA$ 
as well as of $\MD$, which are induced by the filtration given by weight or length, it suffices to understand $ d'(k,l)$.

\begin{prop}\label{prop:core}
We have   for $\MDA$ the identities
\begin{align*}
   \dim_\Q \grw_k  ( \MDA )&= \sum_{i=0}^k d'(k,i) \\
  \dim_\Q \filw_k(\MDA) &= \sum_{j=0}^k \sum_{i=0}^j d'(j,i)\\
  \dim_\Q \filwle_{k,l}(\MDA) &= \sum_{j=0}^k \sum_{i=0}^{l} d'(j,i)
\end{align*}
and for $\MD$ we have
\begin{align*}
\dim_\Q \grwl_{k,l}(\MD)&=  \sum_{j=0}^{k} d'(k-j,l-j)\\ 
 \dim_\Q \grw_k(\MD)&= \dim_\Q \filw_k(\MDA) =\sum_{l=0}^k \sum_{j=0}^{k} d'(k-j,l-j) \\
  \dim_\Q   \filw_k  ( \MD ) &= \sum_{j=0}^k \sum_{i=0}^k \sum_{r=0}^j d'(j-r,i-r) \\
  \dim_\Q   \filwle_{k,l}  ( \MD ) &= \sum_{j=0}^k \sum_{i=0}^l \sum_{r=0}^{j} d'(j-r,i-r)
\end{align*}
\end{prop}
 
\begin{prf}  If $V$ is a vector space with filtration $\operatorname{F}_{\bullet}$ such that
\[
0=\operatorname{F}_0 (V) \subseteq \operatorname{F}_1 (V)\subseteq \dots \subseteq \operatorname{F}_k(V) \subseteq \dots \subseteq V\,,
\]
then $\operatorname{F}_k(V) \cong \oplus_{j\le k} \gr^{\operatorname{F}}_j(V)$. We further 
 know that 
\[\MD \cong \MDA[[1]] \,, \]
 hence modulo $\filwle_{k,l-1}(\MD)$ and  $\filwle_{k-1,l}(\MD)$ we have
 \[
\grwl_{k,l}(\MD) \equiv  \sum_{i=0}^{k} \grwl_{k-i,l-i}(\MDA) [1]^i  . \]
Now the claim follows by the  properties of the product. 
\end{prf}

\begin{thm}
\label{thm-filklmda} We have the following results for $\dim_\Q \filwle_{k,l} (\MDA) $
\begin{small}
\begin{table}[H]\footnotesize
\begin{center}
\begin{tabular} { c | c | c |c|c|c|c|c|c|c|c|c|c|c|c}
$k \backslash l$&0&1&2&3&4&5&6&7&8&9&10&11 \\ \hline
0 &{\bf1}	& 		&			&			& 		& 		& 		&			& 	  & 		& 		&\\ \hline
1 &{\bf1}	&{\bf1}		&			&			& 		& 		& 		& 		& 		& 	  & 		&\\ \hline
2 &{\bf1}	&{\bf2}		&{\bf2}		&			& 		& 		& 		& 		& 		&  		& 		&\\ \hline
3 &{\bf1}	&{\bf3}		&{\bf4}		&{\bf 4}		&& 		&  		& 		&		  & 		& 		&\\ \hline
4 &{\bf1}	&{\bf4}	&{\bf6}		&{\bf 7}&{\bf 7}		& 		& 		& 		&  		& 		& 		&\\ \hline
5 &{\bf1}	&{\bf5}		&{\bf9}		&{\bf12}	&{\bf13}& {\bf13}	& 		& 		& 		& 		& 	&\\ \hline
6 &{\bf1}	&{\bf6}	&{\bf12}	&{\bf18	}	&\blue{22}		& \blue{23}	& \blue{23}	&  & & 		& 		&\\ \hline
7 &{\bf1}	&{\bf7}		&{\bf16}	&\blue{26}	&\blue{35}	& \blue{40}	& \blue{41}	& \blue{41}	&  & &  &\\ \hline
8 &{\bf1}	&{\bf8}	&{\bf20}&\blue{36}	&\blue{53}	& \blue{66}	& \blue{72}	& \blue{73}	& \blue{73}	& & &\\ \hline
9 &{\bf1}	&{\bf9}	&{\bf25}&\blue{48}	&\blue{76}	& \blue{103}& \blue{121}&\blue{128}	& \blue{129}& \blue{129}	& 		&		\\ \hline
10&{\bf1} 	&{\bf10}		&{\bf30}		&\blue{63}		&\blue{107} 	& \blue{155} & \blue{196}	&\blue{220} 	& \blue{228} & \blue{229} & \blue{229} &		\\ \hline
11&{\bf1}	&{\bf11}		&{\bf36}		&\blue{80}		& \blue{145}  	& \blue{225} 	& \blue{304}& \blue{364} 	& \blue{395}	  & \blue{404}   &   \blue{405} &	\blue{405}	 \\ \hline
12&{\bf1}  	&{\bf12}		&{\bf42}		&\blue{100}	& \blue{193}  	& \blue{317} 	& \blue{456}		& ? 	& ?	  & ?   &   ? &	?	 \\ \hline
13&{\bf1}  	&{\bf13}		&{\bf49}		&\blue{123}	& \blue{251}  	& ? 	& ?		& ? 	& ?	  & ?   &   ? &	?	 \\ \hline
14&{\bf1}  	&{\bf14}		&{\bf56}		&\blue{150}	& \blue{321}  	& ? 	& ?		& ? 	& ?	  & ?   &   ? &	?	 \\ \hline
15&{\bf1}  	&{\bf15}		&{\bf64}		&\blue{179}	& ?  	& ? 	& ?		& ? 	& ?	  & ?   &   ? &	?	 
\end{tabular} 
\caption{$\dim_\Q \filwle_{k,l} (\MDA)$:  {\bf proven exact},  {\color{blue} proven lower bounds}.} \label{tab:filklqmz}
\end{center}
\end{table}
\end{small} 
\end{thm}

\begin{prf} 
We first explain how we obtain lower bounds  
 with the help of a computer,
then we give an upper bounds by listing enough relations.

Lower bounds:  
 
We  calculated with the help of a computer  
a reasonable number of the coefficients for each of the brackets in 
$\filwle_{k,l}(\MDA)$.  
 Now the rank of the matrix whose rows are the coefficients gives us  for 
$\dim_\Q \filwle_{k,l}(\MDA)$ a lower bound. 
Since we work only with a finite number of columns, it may happen that we can't distinguish linear independent elements.
The result of our computer calculations is that all the entries\footnote{
 The total running time on a standard PC for each entry 
was less then 24 hours. We point to the fact, that refinements of our code may give some more entries in the table. } in the table of \thmref{thm-filklmda} are lower bounds.

For example in the case of $\filwle_{4,3}(\MDA)$ we checked that the  following matrix     
\[ 
\begin{pmatrix}
1&3&4&7&6&12&8&15\\
\frac{1}{2}&\frac{5}{2}&5&\frac{21}{2}&13&25&25&\frac{85}{2}\\
\frac{1}{6}&\frac{3}{2}&\frac{14}{3}&\frac{73}{6}&21&42&\frac{172}{3}&\frac{195}{2}\\
0&0&1&2&6&7&15&18\\
0&0&1&3&9&15&30&45\\
0&0&\frac{1}{2}&1&4&\frac{9}{2}&\frac{25}{2}&15\\
0&0&0&0&0&1&2&5\\\end{pmatrix}, 
\]  
%
whose rows are the first $8$ coefficients of the $7$ brackets
\[  [2], [3], [4], [2, 1], [2, 2], [3, 1], [2, 1, 1] \]  
has rank $6$. Thus there are at least $7$ (including the constant) linear independent elements in $\filwle_{4,3}(\MDA)$ and therefore $\dim_\Q \filwle_{4,3} (\MDA)\ge 7$.\\\

Upper bounds:

Because of the identity 
\[
\dim_\Q \filwle_{k,l} (\MDA)  = \sum_{i\le k, j \le l} \dim_\Q \grwl_{i,j} (\MDA)
\] 
it suffices to give upper bounds for  $\dim_\Q \grwl_{i,j} (\MDA)$. We use the bounds given 
by $a(k,l)$ minus the number of known relations between the generators.  
There is at least no relations in the generators of $\grwl_{k,1}(\MDA)$, in fact $[k]$ 
is a generator. 
In $\grwl_{k,2}(\MDA)$ we know by \thmref{thm:dim2} that there 
are at least $\left\lfloor \frac{(k-2)}{2} \right\rfloor$ 
relations in between generators.    
In addition we know by \thmref{thm:dim3} the number of relations in length $3$ for the weights $5$ and $6$.  
Now it is easily checked that the lower and upper bounds coincide for the black marked entries
in the table and hence the theorem is proven.  For example in the case of $\filwle_{4,3}(\MDA)$
we have that
\begin{align*}
\dim_\Q \filwle_{4,3} (\MDA)  &\le \sum_{0 \le k\le 4} \sum_{0\le l \le 3} \dim_\Q \grwl_{k,l} (\MDA)\\
&= 1+ \sum_{2 \le k \le 4} (1 -0)  + \sum_{3 \le k \le 4} \Big( \binom{k-2}{1} - \left\lfloor \frac{(k-2)}{2} \right\rfloor \Big)
+  1 -0 \\
& =1+ 3+ 2+1 =7.
\end{align*}

\end{prf}

Unfortunaly there is no direct way to get the dimension of $\grwl_{k,l}(\MDA)$ with the help of a computer. 
However we can deduce the following conditional result.

 \begin{thm} \label{thm-dklmda} 
  i) We have the following
results for $d'(k,l)= \dim_\Q  \grwl_{k,l} (\MDA)  $
\begin{table}[H] 
\footnotesize
\begin{center}
\begin{tabular} { c | c | c |c|c|c|c|c|c|c|c|c|c|c|c}
$k \diagdown l$\!\!&0&1&2&3&4&5&6&7&8&9&10&11 \\ \hline 
0&{\bf1}&{\bf0}&& &&& &&&&&\\ \hline
1&{\bf0}&{\bf0} &&& &&& &&&& \\ \hline
2&{\bf0}&{\bf1}&{\bf0}&&& &&& &&& \\ \hline
3 &{\bf0 }&{\bf1}&{\bf1}&{\bf0} &&& &&& &&\\ \hline
4&{\bf0}&{\bf1}&{\bf1}&{\bf1}&{\bf0} &&& &&&&\\ \hline
5&{\bf0}&{\bf1}&{\bf2}&{\bf2}&{\bf1}&{\bf0} &&&&&&\\ \hline
6&{\bf0}&{\bf1}&{\bf2}&{\bf3}&\blue{3}&\blue{1}&{\bf0} &&&&&\\ \hline
7 &{\bf0}&{\bf1}&{\bf3}&\blue{4}&\blue{5}&\blue{4}&\blue{1}&{\bf0} &&&&\\ \hline
8 &{\bf0}&{\bf1}&{\bf3}&\blue{6}&\blue{8}&\blue{8}&\blue{5}&\blue{1}&{\bf0} &&&\\ \hline
9&{\bf0}&{\bf1}&{\bf4}&\blue{7}&\blue{11}&\blue{14}&\blue{12}&\blue{6}&\blue{1}&{\bf0 }&&\\ \hline
10 &{\bf0}&{\bf1}&{\bf4}&\blue{10}&\blue{16}&\blue{21}&\blue{23}&\blue{17}&\blue{7}&\blue{1}&{\bf0}&\\ \hline
11 &{\bf0}&{\bf1}&{\bf5}&\blue{11}&\blue{21}&\blue{32}&\blue{38}&\blue{36}&\blue{23}&\blue{8}&\blue{1}&{\bf0}\\ \hline
12 &{\bf0}&{\bf1}&{\bf5}&\blue{14}&\blue{28}&\blue{44}&\blue{60}&?&?&\blue{30}&\blue{9}&\blue{1}\\ \hline
13 &{\bf0}&{\bf1}&{\bf6}&\blue{16}&\blue{35}&?&?&?&?&?&\blue{38}&\blue{10}\\ \hline
14 &{\bf0}&{\bf1}&{\bf6}&\blue{20}&\blue{43}&?&?&?&?&?&?&\blue{ 47}\\ \hline
15 &{\bf0}&{\bf1}&{\bf7}&\blue{21}&?&?&?&?&?&?&?&?\\ \hline
\end{tabular} 
\caption{$\dim_\Q \grwl_{k,l}(\MDA)$: {\bf proven}, {\color{blue}conjectured}.} \label{tab:grwl-qmz}
\end{center}
\end{table}

ii) We have the following results for the number of relations in $\dim_\Q \grwl_{k,l}(\MDA)$
\begin{table}[H]\footnotesize
\begin{center}
\begin{tabular} { c | c |c|c|c|c|c|c|c|c|c|c|c|}
$k \diagdown l$\!\!&1&2&3&4&5&6&7&8&9&10 \\ \hline 
1& &&&&&&&&&\\ \hline
2&{\bf 0}&&&&&&&&&\\ \hline
3&{\bf 0}&{\bf 0}&&&&&&&&\\ \hline
4&{\bf 0}&{\bf 1}&{\bf 0}&&&&&&&\\ \hline
5&{\bf 0}&{\bf 1}&{\bf 1}&{\bf 0}&&&&&&\\ \hline
6&{\bf 0}&{\bf 2}&{\bf 3}&{\color{blue} 1}&\blue{0}&&&&&\\ \hline
7&{\bf 0}&{\bf 2}&{\color{blue}6}&{\color{blue}5}&{\color{blue}1}&{\color{blue} 0}&&&&\\ \hline
8&{\bf 0}&{\bf 3}&{\color{blue}9}&{\color{blue}12}&{\color{blue}7}& {\color{blue}1}&{\color{blue} 0}&&&\\ \hline
9&{\bf 0}&{\bf 3}&{\color{blue}14}&{\color{blue}24}&{\color{blue}21}&{\color{blue}9}&{\color{blue}1}&{\color{blue} 0}&&\\ \hline
10&{\bf 0}&{\bf 4}&{\color{blue}18}&{\color{blue}40}&{\color{blue}49}&{\color{blue}33}&{\color{blue}11}&{\color{blue}1}&{\color{blue} 0}&\\ \hline
11&{\bf 0}&{\bf 4}&{\color{blue}25}&{\color{blue}63}&?&?&?&\blue{13}&{\color{blue}1}&{\color{blue} 0}\\ \hline
12 &{\bf 0}&{\bf 5}&{\color{blue}36}&{\color{blue}16}&?&?&?&?&\blue{15}&\blue{1}\\ \hline
13 &{\bf 0}&{\bf 5}&?&?&?&?&?&?&?&\blue{17}\\ \hline
14 &{\bf 0}&{\bf 6}&?&?&?&?&?&?&?&?\\ \hline
15 &{\bf 0}&{\bf 6}&?&?&?&?&?&?&?&?\\ \hline
\end{tabular} 

\caption{Relations in $\dim_\Q \grwl_{k,l}(\MDA)$: {\bf proven},  {\color{blue}conjectured}.} \label{tab:relationsgrqmz} 
\end{center}
\end{table}
 
\end{thm}
  
\begin{prf} i)
If the dimensions of $\filwle_{k,l}(\MDA)$ are given, then    
\begin{align*}
\dim_\Q \grwl_{k,l}(\MDA) =& \dim_\Q \filwle_{k,l}(\MDA) -
\dim_\Q \filwle_{k-1,l} (\MDA)\\
&-  \dim_\Q  \filwle_{k,l-1}(\MDA)  + \dim_\Q  \filwle_{k-1,l-1}(\MDA),
\end{align*}
because  we have
\begin{align*}
\grwl_{k,l}  (\MDA)  \cong & \fille_l \bigg(\filw_k(\MDA)/ \filw_{k-1}(\MDA)\bigg)\\
& \bigg/
\fille_{l-1} \bigg(\filw_k(\MDA)/ \filw_{k-1}(\MDA)\bigg).
 \end{align*}
Now using \thmref{thm-filklmda} we get all the black marked entries in \tabref{tab:grwl-qmz}. For the conjectured entries in \tabref{tab:grwl-qmz} we assumed that all the entries in \tabref{tab:filklqmz} were exact, except for the diagonals for which we guessed the entries for weight bigger then $11$.

ii)  The number of independent relations we found give all the black marked entries in \tabref{tab:relationsgrqmz}, 
since by i) we know that there aren't more.  The conjectured entries in \tabref{tab:relationsgrqmz} equal the difference
of the number of generators $a(k,l)$ of in $\grwl_{k,l}(\MDA)$ minus the corresponding  dimension conjectured in i). 

\end{prf}

\begin{prf}(of \thmref{thm:lowdim_fil_md}) The entries in \tabref{tab:filwl-md} were  calculated from the
values for $d'(k,l)$ given in \thmref{thm-dklmda}  by means of the formula given in \propref{prop:core}. Actually we have double-checked this table with the computer.
\end{prf}

\begin{rem}     
 Of course a lot of the conjectured relations in the table of \thmref{thm-dklmda} can be obtained by using the methods mentioned in this paper. 
We expect that with a more detailed study of the kind of relations 
we can obtain so far we could derive much better results and we plan 
to come back to this in future \cite{Ba:PhDthesis}.   
 

 \end{rem}

\begin{rem} \label{rem:conj-qmz}
The lower bounds where proven with the help of a computer and we expect that our program has found all the
linear independent elements. We therefore conjecture that \tabref{tab:grwl-qmz} in \thmref{thm-dklmda}  gives the exact values 
of $d'(k,l)$ for all $k,l$ we have tested. Assuming this we can ask for relations that are satisfied by
the $d'(k,l)$. We observe that $d'_k = \sum_{l=1}^k d'(k,l)$ satisfies: $d'_0=1$, $d'_1=0$, $d'_2=1$ 
and

\begin{align*} 
d'_k = 2 d'_{k-2} + 2 d'_{k-3} , \quad \textrm{ for } \, 5 \le k \le 11.
\end{align*}
We see no reason why this shouldn't hold for all $k>11$ also, i.e. we ask whether
\begin{align}\label{eq:conj-dk}
\sum_{k\geq 0} \grw_k (\MDA) x^k = \sum_{k\geq 0} d'_k x^k \overset{?}{=} \frac{1-x^2+x^4}{1-2 x^2 - 2 x^3}.
\end{align} 
Even more speculative we may ask whether there a polynomial $P(x,y), Q(x,y) \in \Q[x,y]$ such that
\begin{align}\label{eq:conj-dkl}
\sum_{k , l \geq 0} \dim_Q \grwl_{k,l} (\MDA) x^k y^l = \sum_{k,l\geq 0} d'(k,l) x^k y^l \overset{?}{=} \frac{ P(x,y)}{ Q(x,y)}.
\end{align}
and $\frac{ P(x,1)}{ Q(x,1)} = \frac{1-x^2+x^4}{1-2 x^2 - 2 x^3}$. In fact, for the data we have so far there exist a family of polynomials $ P(x,y)$ and  $Q(x,y)$
such that if $\frac{ P(x,y)}{ Q(x,y)} = \sum a(k,l) x^k y^l$, then 
$d'(k,l)=a(k,l)$ for all $d'(k,l)$ in table in \thmref{tab:grwl-qmz}.

A general reason why such conjectural formulas may hold is that these are analogous to 
the Zagier conjecture for the dimension $d_k$ of $\MZ_k$ 
\[ \sum_{k\geq 0} \dim_Q  \grw_k(\MZ) X^k =
\sum_{k\geq 0} d_k X^k \overset{?}{=} \frac{1}{1-X^2-X^3}
\]
and its refinement by 
the Broadhurst Kreimer conjecture
\[   \sum_{\substack{k\geq 0 \\ l\geq 0}} \dim_\Q    \grwl_{k,l}( \MZ )   X^k Y^l \overset{?}{=} \frac{1 + \mathbb{E}(X) Y}{1 - \mathbb{O}(X) Y + \mathbb{S}(X) Y^2 - \mathbb{S}(X) Y^4} \,. \]
where
\[ \mathbb{E}(X) = \frac{X^2}{1-X^2}  \,,\quad \mathbb{O}(X) = \frac{X^3}{1-X^2} \,, \quad  \mathbb{S}(X) = \frac{X^{12}}{(1-X^4)(1-X^6)} \,.\]

We finally observe that conjecturally the algebra $\MDA$ is much bigger than $\MZ$ as we read of the following table.
\begin{table}[H]\footnotesize
\begin{center}
\begin{tabular} { c| c| c | c | c |c|c|c|c|c|c|c|c|c|c|c|c|c}
$k$ 		&	0 & 1	&	2 & 3 & 4	& 5 & 6 &  7 & 8 & 9 & 10 & 	11 & 12 & 13 & 14 & 15 & 16 \\ \hline
$d_k$  & 1 & 0 & 1 & 1 & 1 & 2 & 2 & 3 & 4 & 5 & 7 & 9 & 12 & 16 & 21 & 28 & 37	\\ \hline
$d'_k$  & 1 & 0 & 1 & 2 & 3 & 6 & 10 & 18 & 32 & 56 & 100 & 176 & 312 & 552 & 976 & 1728 & 3056	\\ 
\end{tabular} 
\caption{First values of $d_k$ and $d'_k$.} 
\end{center}
\end{table}

 \end{rem}

\section{Interpretation as a $q$-analogue of multiple zeta values}

We will show that the brackets can be seen as an q-analogue of multiple zeta values.

\begin{rem} \label{rem:q-vergleich}
The most common example for an q-analogue of multiple zeta values are
the multiple q-zeta values (see for example \cite{db}). 
They are defined for $s_1>1$,$s_2,\dots,s_l \geq 1$ as
\begin{align}\label{eq:defqzeta}
 \zeta_q(s_1,\dots,s_l) = \sum_{n_1 > \dots > n_l > 0} \prod_{j=1}^n \frac{q^{n_j (s_j-1)}}{[n_j]_q^{s_j}}\,, 
\end{align}
where one has to be careful with the notation here, because the brackets $[n]_q$ in this case denote the q-analogue of a natural number $n_j$. They are  given by
\[ [n]_q = \frac{1-q^n}{1-q} = \sum_{j=0}^{n-1} q^j \,. \]
With this it is easy to see that since $s_1>1$  
\[ \lim\limits_{q \rightarrow 1}  \zeta_q(s_1,\dots,s_l) = \zeta(s_1,\dots,s_l)\,.\]
These function also fulfill a lot of relations from which one can deduce relations of MZV due to the limiting process.

It seems  strange to us, that
albeit the cases $(1-q)^s [s]$ have been treated as $q$-zeta values 
\cite{Zud}, \cite{Pupy} or \cite{KKW}
the definition \eqref{eq:defqzeta} has become standard 
(see e.g. \cite{db}, \cite{Zhao:q_mzv},\cite{OKZ})
and not  $(1-q)^{s_1+,\dots,s_k}[s_1,\dots,s_l].$

\end{rem}

\begin{rem}
There is also another $q$-analogue, which is more directly connected to the brackets. It is defined by
\begin{align*}
 \overline{\zeta}_q(s_1,\dots,s_l) &= (1-q)^{-k} \zeta _q(s_1,\dots,s_l) \\
&=  \sum_{n_1 > \dots > n_l >0} \frac{ q^{n_1 (s_1-1)} \dots q^{n_l (s_l-1)}}{(1-q^{n_1})^{s_1}\dots (1-q^{n_l})^{s_l}} 
 \end{align*}
  and which are called modified $q$-multiple zeta values in 
  \cite{OT} or \cite{Tak-harm}.

If all $s_j > 1$, then modified  $q$-multiple zeta values can  be written in terms of brackets, which follows from the fact that the Eulerian polynomials form a basis of a certain space of polynomials 
\cite{BK}. Clearly one has $\overline{\zeta}_q(2,\dots,2) = [2,\dots,2]$ because $P_1(t)=1$. 
If all indices $s_j\ge 2$ the connection gets a little bit more complicated. For example it is
\[  \overline{\zeta}_q(4) = [4] - [3] + \frac{1}{3} [2]  \,, \]
and this is due to the identity
\[ \frac{t^3}{(1-t)^4} =\frac{ t P_3(t) }{3! (1-t)^4 } - \frac{ t P_2(t) }{2! (1-t)^3 } + \frac{1}{3} \frac{ t P_1(t) }{ (1-t)^2 } \,.\]
When one of the $s_j$ is equal to $1$  we don't expect such a simple connection. 
But still there seems to be a connections if $s_1>1$, for example 
\[ \overline{\zeta}_q(2,1) \equiv [2,1] - [2] + \dif[1] \,\mod q^{100}\Q[[q]]\,. \]
It is not difficult to check that the space of modified q-multiple zeta is closed under multiplication
(see e.g. \cite{hi}, p.~2).
However, the algebra of admissible brackets $\MDA$
is not isomorphic to the $\Q$-algebra of modified q-multiple zeta 
values in the sense of   \cite{OT} or \cite{Tak-harm}. This is in essence due to the relation 
$\overline{\zeta}_q(2,1) = \overline{\zeta}_q(3)$ in contrast to $[2,1] \neq [3]$.  
  
\end{rem}

\begin{dfn}\label{def_Zkanal}
For $k \in \N$ we define the map $Z_k : \filw_{k}(\MDA) \rightarrow \R $ by 
\[ Z_k\left( [s_1, \dots , s_l ] \right) = \lim\limits_{q \rightarrow 1}\left((1-q)^k [s_1,\dots,s_l]\right) \,. \]
\end{dfn}

\begin{prop}\label{prop_Zkexpli}
 The map $Z_k$ is linear and on the generators
 of $\filw_{k}(\MDA)$, i.e., on brackets with $s_1>0$, it is given by
\[ Z_k\left( [s_1, \dots , s_l ] \right) = \left\{
\begin{array}{cl} \zeta(s_1,\dots,s_l)\,, & k = s_1+\dots+s_l ,  \\ 0\,, &k > s_1+\dots+s_l 
\,. \end{array}   \right. 
\] 
\end{prop}
\begin{prf} 
Using \lemref{lem:eulerpol} and \lemref{lem:uniformily} below, we derive for $k=s_1+\dots+s_l$
\begin{align*}
 Z_k\left( [s_1, \dots , s_l ] \right) &= \lim\limits_{q \rightarrow 1}\left((1-q)^k [s_1,\dots,s_l]\right) \\
 &=\lim\limits_{q \rightarrow 1}\left((1-q)^k
    \sum_{n_1 > \dots > n_l > 0} \prod_{j=1}^l \frac{q^{n_j} P_{s_j-1}\left( q^{n_j} \right)}{(s_j-1)! (1-q^{n_j})^{s_j}} \right) \\
     &=
    \sum_{n_1 > \dots > n_l > 0} \lim\limits_{q \rightarrow 1} 
     \prod_{j=1}^l \frac{(1-q)^{s_j}}{(1-q^{n_j})^{s_j}} \frac{q^{n_j} P_{s_j-1}(q^{n_j})}{(s_j-1)!  }\\
  &  = \zeta(s_1,\dots,s_l);
 \end{align*}   
here we used that the $k$-th Eulerian polynomial $P_k(t)$ satisfies $P_k(1)=k!$. If $k> s_1 + \dots s_l$ it is $ Z_k\left( [s_1, \dots , s_l ] \right) = \lim\limits_{q \rightarrow 1} (1-q)^{k-s_1-\dots-s_l}  \zeta(s_1,\dots,s_l) = 0$. In \lemref{lem:uniformily} we will justify the interchange of the limit and the summation.
\end{prf}

 \begin{cor}
  Let $f= \sum_{n\geq0} a_n q^n$ be a quasi-modular form of weight $k$. Then the map $Z_k$ sends $f$ to $(-2\pi i)^k a_0$. The space  $S_k$ of weight $k$ cusp-forms is therefore a subspace of the kernel of $Z_k$.\\
\end{cor}
 
\begin{prf} 
 Any quasi-modular form of weight $k$ can be written as a homogenous polynomial in $G_2$, $G_4$ and $G_6$, therefore $\widetilde M_k(\Gamma_1) \subset \MQ_k$. 
Since $Z_k$ is a linear operator we can focus on the monomials. Let us consider the most simplest case first. For $a,b \in \{2,4,6\}$ we have
\[ Z_{a+b}(G_a G_b) = \lim\limits_{q \rightarrow 1} (1-q)^{a+b} G_a G_b = \lim\limits_{q \rightarrow 1} (1-q)^{a} G_a \lim\limits_{q \rightarrow 1} (1-q)^{b} G_b = Z_a(G_a) Z_b(G_b)\, \]
and by \propref{prop_Zkexpli}  we have $Z_a(G_a) Z_b(G_b)=\zeta(a) \zeta(b)$ which is exactly $ (-2 \pi i)^{a+b}$ times the constant term of $G_a G_b$. 
The same argument holds for more general monomials and therefore the claim follows.  
\end{prf}

 \begin{lem} \label{lem:uniformily} \begin{enumerate}[i)]
 \item Define a series $\{F_M(q)\}_{M \in \N}$ by 
 \[ F_M(q) =
 \sum_{M\ge n_1 > \dots > n_l > 0} \prod_{j=1}^l \frac{  (1-q)^{s_j} q^{n_j}P_{s_j-1}\left( q^{n_j} \right) }{(1-q^{n_j})^{s_j} (s_j-1)! }
 \, ,\]
 then  it converges uniformly to $(1-q)^k [s_1,\dots,s_l]$  for $q$ in the interval $[0,1]$ 
 and therefore
 \begin{align}\label{eq:uniform?} \lim_{q \to 1}\left(
(1-q)^k [s_1,\dots,s_l] \right)= \sum_{n_1 > \dots > n_l > 0} \lim_{q \to 1}
 \prod_{j=1}^l \frac{  (1-q)^{s_j} q^{n_j}P_{s_j-1}\left( q^{n_j} \right) }{(1-q^{n_j})^{s_j} (s_j-1)! }
\end{align}
\item 
Let $k,n \in \N$ be natural numbers and define the function
\[ f_{k,n}(q) = \frac{(1-q)^k  q^n P_{k-1}(q^n)}{(1-q^n)^k (k-1)!} \,,\]
then for $q \in [0,1]$ it is
  $f_{1,n}(q) \leq \frac{1}{n}$ 
 and for $k>1$ we have $f_{k,n}(q) \leq \frac{1}{n^2} $.
\end{enumerate}
\end{lem}
\begin{prf} We start with the proof of ii) because we need it for the proof of i). It is
\[ f_{1,n}(q) = \frac{(1-q) q^n}{(1-q^n)} \] 
because $P_0(q^n)=1$. This is bounded by $\frac{1}{n}$ because the function 
\[ b_n(q) = n (1-q)q^n-(1-q^n) \]
is negative for all $q\in (0,1)$ which can be seen by $b_n(1) = 0$ and the fact that the derivative 
\[ b'_n(q) = n^2(1-q)q^{n-1} + n (q^{n-1} - q^n) \,. \]
is positive.
We will show that \[ \frac{(1-q)^2 q^n}{(1-q^n)^2} \cdot \frac{P_{k-1}(q^n)}{(k-1)!} \leq \frac{1}{n^2} \]
for all $k$. This will be sufficient ii) for proving the statement for all $k\geq 2$ because it is $\frac{1-q}{1-q^n}<1$ for $q \in (0,1)$. Because of the positivity of the coefficients of $P_{k-1}(q)$ and $P_{k-1}(1) = (k-1)!$ we have for $q \in (0,1)$ that 
\[ \frac{P_{k-1}(q^n)}{(k-1)!}  \leq 1 \,. \]
It therefore remains to show that 
\[ h_n(q) := \frac{(1-q)^2 q^n}{(1-q^n)^2} - \frac{1}{n^2} \stackrel{!}{\leq} 0 \,. \]
We will do this by showing that $h_n(q)$ is monotonically increasing in the desired interval and 
\[ \lim\limits_{q \rightarrow 1}{ h_n(q) } = 0 \,. \]
The latter can be seen by using l'hospital twice. For the monotonicity we first derive the derivative of $h$: 
\[ h'_n(q) = \frac{-(1-q) q^{n-1}}{(1-q^n)^3} \cdot \big(  2 q (1-q^n) - n (1-q) (1+q^n) \big)   \,. \]
The first factor is negative and we therefore just have to proof that the term in the brackets is also negative for all $n \in \N$ and $q\in (0,1)$ which we will do by induction on n. For $n=1$ this is trivial and for the inductive step we first rewrite the statement as
\[ 2 \frac{q (1-q^n)}{(1-q)} = 2 \sum_{j=1}^n q^j \leq n (1+q^n) \,. \]
Assuming that this holds for an $n$ we can write
\begin{align*}
 2 \sum_{j=1}^{n+1} q^j &=  2 \sum_{j=1}^{n} q^j + 2 q^{n+1} \leq n (1+q^n) + 2 q^{n+1} \,.
\end{align*}
Now we have to show that 
\[ n (1+q^n) + 2 q^{n+1} \stackrel{!}{\leq}  (n+1) (1+q^{n+1}) \]
which we again do by first setting
\[ g_n(q) := (n+1) (1+q^{n+1}) - \left(  n (1+q^n) + 2 q^{n+1} \right) = n (q^{n-1}-q^n) + 1 - q^{n+1}  \]
and then noticing that $g_n(1)=0$. The derivative $g'_n(q) = -q^{n-1}(n^2 (1-q) + q)$ is clearly negative for $q \in (0,1)$ which implies $g_n(q) \geq 0$ and therefore finishes the inductive step.  
 
We now prove i). 
  Using the  bounds in ii) and taking into account $s_1>1$ we have the bound
\begin{align*}
 F_M(q) &=
 \sum_{M\ge n_1 > \dots > n_l > 0} 
\prod_{j=1}^l \frac{  (1-q)^{s_j} q^{n_j}P_{s_j-1}\left( q^{n_j} \right) }{(1-q^{n_j})^{s_j} (s_j-1)! } \\
&\leq 
 \sum_{M\ge n_1 > \dots > n_l > 0} \frac{1}{n_1^{2} n_2 \dots n_l} 
  \leq \zeta(2,1,\dots,1)
  \end{align*}
for $q \in [0,1]$ and all $ M >0$.
Therefore the sum on the right-hand side 
of \eqref{eq:uniform?} converges uniformly as a function in $q$  and therefore we can interchange limit and summation. \end{prf}
 
 
\begin{rem} \label{rem:pupy} In \cite{Pupy} it is shown $[1] \approx - \frac{\log(1-q)}{1-q}$ near $q=1$. Since
a bracket $[s_1,\dots,s_l]$ with $s_1=1$ are polynomials in $[1]$, it is clear that 
$Z_k$ can't be extended as an analytical map as given in 
\defref{def_Zkanal} to all $\filw_{k}(\MD)$.
\end{rem}

\section{Applications to multiple zeta values}

As mentioned in the introduction we now want to consider a direct connection of brackets with multiple zeta values (MZV). 

We start by defining    for any $\rho \in \R_{\geq 1}$ the following spaces
\[ \MQ_\rho = \left\{ \sum_{n>0} a_n q^n \in \R[[q]] \,  \mid \, a_n = O(n^{\rho-1}) \right\}  \]
and
\[ \MQ_{<\rho} = \left\{ \sum_{n>0} a_n q^n \in \R[[q]] \, \mid \, \exists\, \varepsilon>0 \text{ with } a_n = O(n^{\rho-1-\varepsilon}) \right\} \,, \]
 where $a_n = O(n^{\rho-1})$ is the usual big $O$ notation which means that there is an $C \in \R$ with $| a_n | \leq C n^{\rho-1}$ for all $n \in \N$.

\begin{lem} \label{lem_Qrho}
\begin{enumerate}[i)]
\item Both  $\MQ_{<\rho}$ and $\MQ_\rho$ are $\R$ vector spaces.
\item We have $\MQ_{\rho-1} \subset \MQ_{<\rho} \subset \MQ_\rho$. 
\item Let $r,s \in \R_{\geq 1}$ then 
\[ \MQ_{<r} \cdot \MQ_{<s} \subset \MQ_{<r+s},\,\, \MQ_{<r} \cdot \MQ_s \subset \MQ_{<r+s} 
\textrm{ and } \MQ_r \cdot \MQ_s \subset \MQ_{r+s}\,. \]
\end{enumerate}
\end{lem}
\begin{prf} It is obvious that i) and ii) hold. For iii) we consider
 $f = \sum_{n>0} a_n q^n \in \MQ_r$, $g=\sum_{n>0} b_n q^n \in \MQ_s$. Then by definition $|a_n| \leq C_1 n^{r-1}$ and $|b_n| \leq C_2 n^{s-1}$ for some constants $C_1$ and $C_2$. Setting $f \cdot g = \sum_{n>0} c_n q^n$ we derive
\[ |c_n| = \left| \sum_{n_1+n_2 = n} a_{n_1} b_{n_2} \right| \leq C_1 C_2  \sum_{n_1+n_2 = n} n_1^{r-1} n_2^{s-1} \leq C_1 C_2 n \cdot n^{r-1} n^{s-1} = O(n^{r+s-1})  \,. \]
and therefore $f \cdot g \in \MQ_{r+s}$.  By similar considerations the remaining cases follow.
\end{prf}

\begin{prop} \label{pro_Qrho}For $\rho>1$
define the map $Z_\rho$ for a  $f = \sum_{n>0} a_n q^n  \in \R[[q]]$ by 
\[ Z_\rho(f) = \limsup\limits_{q \rightarrow 1}{(1-q)^\rho  \sum_{n>0} a_n q^n} \,, \]
where one assumes $q \in (0,1)$.
Then the following statements are true
\begin{enumerate}[i)]
\item $Z_\rho$ is a linear map from $\MQ_\rho$ to $\R$ 
\item $\MQ_{<\rho}  \subset \ker{Z_\rho}$.  
\item $\dif \MQ_{<\rho-1} \subset  \ker(Z_{\rho})$, where as before $\dif= q \frac{d}{dq}$.
\end{enumerate}
\end{prop}

\begin{prf} We prove i) and ii) simultaneously.  
 In order to do this we use the following expression for the polylogarithm 
\[ \Li_{-s}(q) = \Gamma(1+s) (-\log q)^{-s-1} + \sum_{n=0}^\infty \frac{\zeta(-s-n)}{n!} \,(\log q)^n  \]
which is valid for $s \neq -1,-2,-3,\dots$, $|z|<1$ and where $\zeta(-s-n)$ is the analytic continuation of the Riemann zeta-function. The proof of this can be found in \cite{og} Corollary 2.1. 
The logarithm has the following expansion near $q=1$
\[ - \log(q) = \sum_{n=1}^\infty \frac{(1-q)^n}{n} \,.\]
Using this one gets for $\varepsilon\geq 0$
\begin{align*}
 &\limsup\limits_{q \rightarrow 1}{(1-q)^\rho  \sum_{n>0} n^{\rho-1-\varepsilon} q^n}  =   \lim\limits_{q \rightarrow 1}{(1-q)^\rho  \Li_{\varepsilon + 1 - \rho}(q)} \\
&=   \limsup\limits_{q \rightarrow 1}{(1-q)^\rho  \left( \Gamma(\rho-\varepsilon) (-\log q)^{-\rho+\varepsilon} + \sum_{n=0}^\infty \frac{\zeta(-\rho+\varepsilon-n)}{n!} \, (\log q)^n \right) } \\
&= \Gamma(\rho-\varepsilon)  \limsup\limits_{q \rightarrow 1}{ \frac{(1-q)^\rho}{ \left( \sum_{n=1}^\infty \frac{(1-q)^n}{n} \right)^{\rho-\varepsilon} } } = \left\{\begin{array}{cl} \Gamma(\rho) \,, & \varepsilon = 0 \\ 0\,, & \varepsilon>0 \end{array}\right.  \,.
\end{align*}
Now assume that for a $\varepsilon \geq 0$ we have $f = \sum_{n>0} a_n q^n$ with $|a_n|\leq C \cdot n^{\rho-1-\varepsilon}$, i.e. $f \in \MQ_\rho$ for $\varepsilon = 0$ and  $f \in \MQ_{<\rho}$ for $\varepsilon > 0$, then the calculation above gives
\[ | Z_\rho(f)| = \left|  Z_\rho \left(  \sum_{n>0} a_n q^n  \right) \right| \leq C \cdot Z_\rho \left( \sum_{n>0} n^{\rho-1-\varepsilon} \right) = \left\{\begin{array}{cl} C \cdot \Gamma(\rho) \,, & \varepsilon = 0 \\ 0\,, & \varepsilon>0 \end{array}\right. \] 
and therefore $Z_\rho(f) \in \R$ and $Z_\rho(f) = 0$ respectively. 
 
For iii) we just have to observe that the derivative $\dif = q \frac{d}{dq}$ on $\sum_{n>0} a_n q^n$ is given by $\sum_{n>0} n a_n q^n$. With this it is clear that with i) we obtain $\dif \left( \Q_{< \rho-1} \right) \subset \Q_{<\rho} \subset \ker(Z_\rho)$.  
\end{prf}

The brackets $[s_1,\dots,s_l]$ can be considered as elements in the spaces we studied above.  

\begin{prop}\label{prop_Zk-theo}
\begin{enumerate}[i)]
\item For any $s_1,\dots,s_l$ we have $[s_1,\dots,s_l] \in \MQ_{< s_1 + \dots + s_l + 1}$.
\item If all $s_1,s_2,\dots,s_l > 1$, then $[s_1,\dots,s_l] \in \MQ_{s_1 + \dots + s_l}$.
\item For any $s_1,\dots,s_l$ we have  \[
 [s_1,\dots,s_l]   \subset  \ker(Z_{s_1+\dots+s_l+1})\]
and 
\[\dif [s_1,\dots,s_l]   \subset  \ker(Z_{s_1+\dots+s_l+2}) \]
\end{enumerate}
\end{prop}

\begin{prf} We begin with the proof of ii). 
It is a well-known fact that for $s>1$ the divisor sums $\sigma_{s-1}(n)$ are in $O(n^{s-1})$ and therefore $[s] \in \MQ_s$. Then by \lemref{lem_Qrho} iii)
we have $\sum_{n>0} a_{s_1,\dots,s_l}(n) q^n := [s_1]\dots [s_l] \in \MQ_{s_1+\dots+s_l}$.
  
It is clearly 
\begin{align*}
\sigma_{s_1-1,\dots,s_l-1}(n)  &= \sum_{\substack{u_1 v_1 + \dots + u_l v_l = n\\u_1 > \dots > u_l >0}} v_1^{s_1-1} \dots v_l^{s_l-1}\\
& \leq \sum_{u_1 v_1 + \dots + u_l v_l = n} v_1^{s_1-1} \dots v_l^{s_l-1}   = a_{s_1,\dots,s_l}(n) \,, 
\end{align*}
which implies $[s_1,\dots,s_l] \in \MQ_{s_1+\dots+s_l}$.
 
In order to show i) we can use the same argument as in ii) except that   one has $\sigma_0(n) \in O(\log(n)) \subset O(n^\varepsilon)$ for any $\varepsilon>0$. Using this we obtain $[1] \in \Q_{<2}$ and therefore $[s_1,\dots,s_l] \in \Q_{<s_1+\dots+s_l+1}$ for $s_1,\dots,s_l \geq 1$.

Finally iii) is an immediate consequence of \propref{pro_Qrho} ii) and iii). 
\end{prf}

Using $\MD = \MDA[[1]]$ we define a map 
\begin{align*}
Z^{alg}_k&:  \filw_{k}(\MD) \longrightarrow \R[T] \,,\\ 
Z^{alg}_k&\left( \sum_{j=0}^{k} g_j [1]^{k-j} \right) = \sum_{j=0}^{k} Z_j(g_j) T^{k-j} \in \R[T]\, ,
\end{align*}
where $g_j \in \filw_{j}(\MDA)$. 

\begin{prop} \label{prop:zkzkalg}
For all $f \in \filw_{k}(\MD)$ it is $Z_{k+2}^{alg}(\dif f) = 0$. 
\end{prop}
\begin{prf}
An element in  $\filw_{k}(\MD)$ can be written as $\sum_{j=0}^{k} g_j [1]^{k-j} $ with $g_j \in \filw_{j}(\MDA)$. The map $Z^{alg}_k$ is linear,
it therefore suffices to prove the statement for a $f \in \filw_{k}(\MD)$ of the form $f= g_j [1]^{k-j}$. The derivative of this $f$ is given by
\[ \dif f = \dif g_j \cdot [1]^{k-j} + (k-j) g_j \dif[1]\cdot  [1]^{k-j-1} \,.\] 
As we saw before it is $\dif [1] = [3] + \frac{1}{2}[2] - [2,1] \in \MDA$ and by \propref{prop:qmzcloseddif} it is $\dif g_j \in \MDA$. The map $Z_{k+2}^{alg}$ is therefore given on $\dif f$ by
\[ Z_{k+2}^{alg}(\dif f) = Z_{j+2}\left( \dif g_j \right) T^{k-j} + (k+j) Z_{j+3}\left( g_j \dif[1] \right) T^{k-j-1} \,.  \] 
It is $ Z_{j+3}\left( g_j \dif[1] \right) = Z_j(g_j) \cdot Z_3(\dif [1])$ and by \propref{prop_Zk-theo} we obtain $Z_3\left(\dif[1] \right) = Z_{j+2}\left( \dif g_j \right) = 0$ from which the statement follows
\end{prf}

\begin{rem}
The authors also expect that  the implication 
\[
Z_k(f) = 0 \Longrightarrow Z^{alg}_k(f) = 0
\]
holds for arbitrary $f \in \MD$. 
\end{rem}

Now \thmref{thm_Zk-kernel} follows by \propref{prop_Zk-theo} and \propref{prop:zkzkalg}. Using these propositions we are able to derive relations between MZV coming from elements in the kernel of the map $Z_k$. We give a few examples which give a new interpretation of well known identities of multiple zeta values.

\begin{ex} 
\begin{enumerate}[i)]
\item We have seen earlier that the derivative of $[1]$ is given by
\[ \dif [1] = [3] + \frac{1}{2}[2] - [2,1] \]
and because of the proposition it is $\dif[1], [2] \in \ker Z_3$ from which $\zeta(2,1) = \zeta(3)$ follows. 
\item (Shuffle product) \propref{cor:formularfordk} stated for $s_1 + s_2 = k+2$ that
\[ \binom{k}{s_1-1} \frac{d[k]}{k} = [s_1]\cdot [s_2] +\binom{k}{s_1-1} [k+1] - \sum_{a+b=k+2} \left( \binom{a-1}{s_1-1}+\binom{a-1}{s_2-1} \right) [a,b] \,.\]
Applying $Z_{k+2}$ on both sides we obtain the shuffle product for single zeta values
\[ \zeta(s_1) \cdot \zeta(s_2) = \sum_{a+b=k+2} \left( \binom{a-1}{s_1-1}+\binom{a-1}{s_2-1} \right) \zeta(a,b) \,.\]
\end{enumerate}
\end{ex}

\begin{ex}
For the cusp form $\Delta \in S_{12} \subset \ker(Z_{12})$ we derived the representation 
\begin{align*}
\frac{1}{2^6 \cdot 5 \cdot 691} \Delta  &=  168 [5,7]+150 [7,5]+28 [9,3]\\
&+\frac{1}{1408} [2] - \frac{83}{14400}[4] +\frac{187}{6048} [6] - \frac{7}{120} [8] - \frac{5197}{691} [12] \,.
\end{align*}
Letting $Z_{12}$ act on both sides one obtains the relation
\[ \frac{5197}{691} \zeta(12) =  168 \zeta(5,7)+150 \zeta(7,5) + 28 \zeta(9,3) \,. \] 
In general it is known due to \cite{gkz} that every cusp form of weight $k$ give rise to a relation between double zeta values with odd entries modulo $\zeta(k)$. We believe that one can give an alternative proof of his fact with the help of brackets.
\end{ex}

At the end we want to mention a curious property of the brackets which seems to appear at length $3$. Fixing a weight $k$ and a length $l$ one could ask, if there are linear relations 
between brackets $[s_1,\dots,s_l]$ with the \underline{same} weight $s_1+\dots+s_l = k$ and length $l$.
For $l=2$ using the computer the authors could not find any such relations up to weight $30$. But for $l=3$ there seem to be relations of this form starting in weight $9$. The first two of them are given by:
\begin{conj}
In $\filwle_{9,3}(\MD)$ and $\filwle_{10,3}(\MD)$ we have the relation
\begin{align*}
0 \,=&\, \frac{9}{5}\,[2, 3, 4]+2\,[2, 4, 3]-[2, 5, 2]\\
&+2 \left( [3, 5, 1]-[3, 1, 5] \right) - \frac{1}{5}\,[3, 2, 4]-[3, 3, 3]-[3, 4, 2]\\
&+\frac{3}{5} \left( [4, 4, 1]-[4, 1, 4] \right) -\frac{11}{10}\,[4, 2, 3]+\frac{1}{2}\,[4, 3, 2]\\
&+\frac{4}{5} \left( \,[5, 1, 3]-[5, 3, 1] \right)-[6, 1, 2]+[6, 2, 1]\,.
\end{align*}
and 
\begin{align*}
0 \,=&\, \frac{4}{3}\,[2, 3, 5]+\frac{14}{5}\,[2, 4, 4]+\frac{29}{15}\,[2, 5, 3]-[2, 6, 2]\\
&+2 \left( [3, 6, 1]- [3, 1, 6]\right) - \frac{2}{3} \,[3, 2, 5]+\frac{2}{5} \,[3, 3, 4]-\frac{1}{15} \,[3, 4, 3]-[3, 5, 2]\\
&+2 \left( [4, 5, 1]- [4, 1, 5]\right) - \frac{6}{5}\,[4, 2, 4]-\frac{4}{3}\,[4, 3, 3]-\frac{2}{5}\,[4, 4, 2]\\
&+\frac{2}{5} \left( [5, 4, 1]- [5, 1, 4] \right)-[5, 2, 3]+\frac{1}{5}\,[5, 3, 2]\\
&+\frac{1}{3} \left([6, 1, 3]- [6, 3, 1] \right) -[7, 1, 2]+[7, 2, 1] \,.
\end{align*}
\end{conj}
Notice that these are all elements in $\filw_{9} (\MDA)$ (resp. $\filw_{10} (\MDA)$)  and therefore a relation for triple zeta values would follow from this. There are similar relations in higher weights and computations show the following: 
 \begin{table}[H]
\begin{center}
\begin{tabular} {c|c|c|c|c|c|c|c|c|c|c|c|c|c}
$k$   & 1-8 & 9 & 10 & 11 & 12   & 13  & 14   &  15   & 16  & 17  & 18  &  19  & 20 \\ \hline
$t_k$ & 0   & 1 & 1  &  3  & 6   & 8   &  12   & 16   & 21  & 25  & 32  &  37	 &  45   \\
\end{tabular} 
\caption{Conjectured numbers $t_k$ of relations between $[a,b,c]$ with $a+b+c = k$.} 
\end{center}
\end{table}

{\small
{\it E-mail:}\\\texttt{henrik.bachmann@math.uni-hamburg.de}\\ \texttt{kuehn@math.uni-hamburg.de}\\

\noindent {\sc Fachbereich Mathematik (AZ)\\ Universit\"at Hamburg\\ Bundesstrasse 55\\ D-20146 Hamburg}}

\begin{thebibliography}{40}

\bibitem[AR]{ar}  G. Andrews, S. Rose:
{\it MacMahon's sum-of-divisors functions, Chebyshev polynomials, and Quasi-modular forms}, preprint, arXiv:1010.5769 [math.NT].

\bibitem[Ba1]{hb} H. Bachmann:
{\it Multiple Zeta-werte und die Verbindung zu Modulformen durch Multiple Eisensteinreihen}, Master thesis, Universit\"at Hamburg 2012. \url{http://www.math.uni-hamburg.de/home/bachmann/msc_henrik_bachmann.pdf}.

\bibitem[Ba2]{Ba:PhDthesis} H. Bachmann: 
{\it PhD-Thesis}, in preparation, Universit\"at Hamburg.

\bibitem[BBK]{BBK} H. Bachmann, O. Bouillot, U. K\"uhn:
{\it The algebra of multiple divisor functions and applications to multiple Eisenstein series}, work in progress.

\bibitem[BT]{BT} H. Bachmann, K. Tasaka:
{\it Multiple Eisenstein series and the Goncharov coproduct}, work in progress.

\bibitem[BK]{BK} H. Bachmann, U. K\"uhn:
{\it  A short note on a conjecture of Okounkov about a $q$-analogue of multiple zeta values}, work in progress.

\bibitem[Br]{db} D. Bradley:
{\it Multiple q-zeta values}, J. Algebra 283 (2005), no. 2, 752-798. 

\bibitem[CG]{og}  O. Costin, S. Garoufalidis:
{\it Resurgence of the fractional polylogarithms}, Mathematical Research Letters 16 (2009), Nr. 5, 817-826.
 
\bibitem[Fo]{df} D. Foata:
{\it Eulerian polynomials: from Euler's Time to the Present}, The legacy of Alladi Ramakrishnan in the mathematical sciences, 253-273, Springer, New York, 2010. 

\bibitem[GKZ]{gkz} H. Gangl, M. Kaneko, D. Zagier:
{\it Double zeta values and modular forms}, in "Automorphic forms and zeta functions" p. 71-106, World Sci. Publ., Hackensack, NJ, 2006.

\bibitem[H1]{h1} M. Hoffman:
{\it The Algebra of Multiple Harmonic Series}, Journal of Algebra, Volume 194, Issue 2, 15 August 1997, Pages 477-495, ISSN 0021-8693, 10.1006/jabr.1997.7127.

\bibitem[HI]{hi} M. Hoffman, K. Ihara:
{\it Quasi-shuffle products revisited}, Max-Planck-Institut f�r Mathematik Preprint Series 2012 (16).

\bibitem[HO]{HofOhn} M. Hoffman, Y. Ohno: 
{\it Relations of multiple zeta values and their algebraic expression}, J. Algebra 262(2), 332?347 (2003).

\bibitem[IKZ]{ikz}  K. Ihara, M. Kaneko, and D. Zagier:
{\it Derivation and double shuffle relations for multiple zeta values}, Compositio Math. 142 (2006), 307-338.

\bibitem[KKW]{KKW} M. Kaneko, N. Kurokawa, M. Wakayama: 
{\it A variation of Euler's approach to values of the Riemann zeta function}, Kyushu J. Math., 57 (2003), 175-192.

\bibitem[Ma]{mm} P. A. MacMahon:
{\it Divisors of numbers and their continuations in the theory of partitions}, Reprinted: Percy A. MacMahon Collected Papers (G. Andrews, ed.), MIT Press, Cambridge, 1986, pp. 305-341.

\bibitem[OKZ]{OKZ} Y. Ohno, J. Okuda, and W. Zudilin:
{\it Cyclic q-MZSV sum}, J. Number Theory 132 (2012), 144-155.

\bibitem[OT]{OT} J. Okuda, Y. Takeyama:
{\it On relations for the multiple q-zeta values}, Ramanujan J. 14 (2007), 379-387.

\bibitem[OZ]{ZagOhn} Y. Ohno, D. Zagier: 
{\it Multiple zeta values of fixed weight, depth, and height}, Indag. Math. (N.S.) 12(4), 483-487 (2001).

\bibitem[Pu]{Pupy} Yu. A. Pupyrev:
{\it Linear and Algebraic Independence of q-Zeta Values},Mathematical Notes, vol. 78, no. 4, 2005, pp. 563-568.
Translated from Matematicheskie Zametki, vol. 78, no. 4, 2005, pp. 608-613.

\bibitem[Ta]{Tak-harm} Y. Takeyama:
{\it The algebra of a q-analogue of multiple harmonic series}, preprint	arXiv:1306.6164 [math.NT].

\bibitem[Zh]{Zhao:q_mzv} J. Zhao:
{\it  Multiple q-zeta functions and multiple q-polylogarithms}, Ramanujan J. 14 (2007), 189-221.
 
\bibitem[Za1]{zag77} D. Zagier:
{\it Modular forms whose Fourier coefficients involve zeta-functions of quadratic fields}, In "Modular functions of one variable", VI (Proc. Second Internat. Conf., Univ. Bonn, Bonn, 1976), pages 105-169. Lecture Notes in Math., Vol. 627. Springer, Berlin, 1977.

\bibitem[Za2]{dz} D. Zagier:
{\it Elliptic modular forms and their applications}, In "The 1-2-3 of modular forms", 1-103, 
Universitext, Springer, Berlin, 2008. 


\bibitem[Zu]{Zud} V. V. Zudilin:
{\it Diophantine Problems for q-Zeta Values}, Mathematical Notes, vol. 72, no. 6, 2002, pp. 858-862.
Translated from Matematicheskie Zametki, vol. 72, no. 6, 2002, pp. 936-940.


\end{thebibliography}
\end{document}